\RequirePackage{ifpdf}
\ifpdf 
\documentclass[pdftex]{sigma}
\else
\documentclass{sigma}
\fi

\numberwithin{equation}{section}

\newtheorem{Theorem}{Theorem}[section]
\newtheorem*{Theorem*}{Theorem}
\newtheorem{Corollary}[Theorem]{Corollary}
\newtheorem{Lemma}[Theorem]{Lemma}
\newtheorem{Proposition}[Theorem]{Proposition}
{ \theoremstyle{definition}
	\newtheorem{Definition}[Theorem]{Definition}
	
	\newtheorem{Example}[Theorem]{Example}
	\newtheorem{Remark}[Theorem]{Remark} }

\usepackage{tikz}
\usepackage{tikz-cd}
\usepackage[shortlabels]{enumitem}

\usetikzlibrary{decorations.pathreplacing,calligraphy}

\DeclareMathOperator{\Hom}{{Hom}}

\DeclareMathOperator{\Newt}{{Newt \,}}

\DeclareMathOperator{\Pic}{{Pic}}

\DeclareMathOperator{\Supp}{{Supp}}

\DeclareMathOperator{\cD}{\mathcal{D}}

\DeclareMathOperator{\cO}{\mathcal{O}}

\DeclareMathOperator{\cX}{\mathcal{X}}
\DeclareMathOperator{\cY}{\mathcal{Y}}

\newcommand{\CC}{\mathbb{C}}
\newcommand{\FF}{\mathbb{F}}

\newcommand{\PP}{\mathbb{P}}
\newcommand{\QQ}{\mathbb{Q}}
\newcommand{\RR}{\mathbb{R}}

\newcommand{\ZZ}{\mathbb{Z}}

\newcommand{\vgap}{\vspace{0.1em}}

\makeatletter
\newcommand{\xleftrightarrow}[2][]{\ext@arrow 3359\leftrightarrowfill@{#1}{#2}}
\newcommand{\xdashrightarrow}[2][]{\ext@arrow 0359\rightarrowfill@@{#1}{#2}}
\newcommand{\xdashleftarrow}[2][]{\ext@arrow 3095\leftarrowfill@@{#1}{#2}}
\newcommand{\xdashleftrightarrow}[2][]{\ext@arrow 3359\leftrightarrowfill@@{#1}{#2}}
\def\rightarrowfill@@{\arrowfill@@\relax\relbar\rightarrow}
\def\leftarrowfill@@{\arrowfill@@\leftarrow\relbar\relax}
\def\leftrightarrowfill@@{\arrowfill@@\leftarrow\relbar\rightarrow}
\def\arrowfill@@#1#2#3#4{%
	$\m@th\thickmuskip0mu\medmuskip\thickmuskip\thinmuskip\thickmuskip
	\relax#4#1
	\xleaders\hbox{$#4#2$}\hfill
	#3$%
}
\makeatother

\begin{document}
	
\allowdisplaybreaks

\newcommand{\arXivNumber}{2112.08246}

\renewcommand{\PaperNumber}{095}

\FirstPageHeading
	
\ShortArticleName{Mirrors to Del Pezzo Surfaces and the Classification of $T$-Polygons}
	
\ArticleName{Mirrors to Del Pezzo Surfaces\\ and the Classification of $\boldsymbol{T}$-Polygons}
	
\Author{Wendelin LUTZ}
	
\AuthorNameForHeading{W.~Lutz}
	
\Address{Department of Mathematics and Statistics, Lederle Graduate Research Tower,\\
		University of Massachusetts, Amherst, MA 01003-9305, USA}
	\Email{\href{mailto:wendelinlutz@umass.edu}{wendelinlutz@umass.edu}}
	\URLaddress{\url{https://sites.google.com/view/wendelin-lutz}}
	
\ArticleDates{Received May 07, 2024, in final form October 14, 2024; Published online October 22, 2024}
	
\Abstract{We give a new geometric proof of the classification of $T$-polygons, a theorem originally due to Kasprzyk, Nill and Prince, using ideas from mirror symmetry. In particular, this gives a completely geometric proof that any two toric $\mathbb{Q}$-Gorenstein degenerations of a~smooth del Pezzo $X$ surface are connected via trees of rational curves in the moduli space~of~$X$.}
	
	
\Classification{14J33; 14E07}

\section{Introduction}
Let $X$ be a Fano manifold of dimension $n$. It is an open problem to classify the possible deformation types of Fano manifolds, and as a first step, it is natural to consider degenerations of $X$ to a toric Fano variety. More precisely, following \cite{AkhtarCoatesCorti}, we consider normal degenerations of $X$ to a toric Fano variety such that $-K$ is relatively ample and $\QQ$-Cartier, usually called $\QQ$-Gorenstein ($qG$) degenerations.
In the moduli space of $X$, many of these toric degenerations are connected by trees of rational curves: following Ilten \cite{Ilten}, we consider pencils $f \colon \cX \rightarrow \PP^1$ with $f^*(0)$ and $f^*(\infty)$ toric Fano varieties. In the context of toric degenerations, these pencils arise naturally from a combinatorial operation on polytopes called mutation: let $P$ be a lattice polytope with the origin in its interior, and write $X_P$ for the toric Fano variety associated to the fan over the faces of $P$.\footnote{The moment polytope of $X_P$ is the dual polytope of $P$.} A mutation of $P$ produces another lattice polytope $Q$ with the property that the toric varieties $X_P$ and $X_Q$ are related by a $qG$-pencil $f \colon \cX \rightarrow \PP^1$ with $f^*(0)=X_P$ and~${f^*(\infty)=X_Q}$.
This shows that, if $X$ admits a $qG$-degeneration to $X_P$, we obtain many other $qG$-degenerations $X_Q$ of $X$ which can be connected by trees of rational curves via these special $qG$-pencils.

One might hope that all toric degenerations of $X$ are connected in this way. If true, this would give a bijection between mutation equivalence classes of lattice polytopes and deformation types of those Fano manifolds which admit a $qG$-degeneration to a toric Fano variety.

The evidence in dimension $n=2$ is encouraging: lattice polygons that arise from $qG$-degenerations of a smooth del Pezzo surface $X$ are called $T$-polygons, and the bijection has been verified by Kasprzyk--Nill--Prince \cite{KasprzykNillPrince}: the authors use a combinatorial argument to show that there are precisely ten mutation-equivalence classes of $T$-polygons and these ten classes biject with the ten deformation families of smooth del Pezzo surfaces. In particular, this implies Conjecture A of \cite{AkhtarCoatesCorti} for smooth del Pezzo surfaces:
\begin{Theorem}[Kasprzyk--Nill--Prince]\label{thm:main}
		Let $X$ be a smooth del Pezzo surface, and let $X_P$ and~$X_Q$ be two toric qG-degenerations of $X$. Then $X_P$ and $X_Q$ are connected by a chain of $\PP^1$'s in the moduli space of $X$. More precisely, there exist
		qG-families $f_i \colon \cX_i \rightarrow \PP^1$, $1 \leq i \leq n$,
		such that we have the following equalities of scheme-theoretic inverse images
		\[
		f^*_1(0)=X_P, \qquad f^*_i(\infty)=f_{i+1}^*(0), \qquad f_n^*(\infty)=X_Q.
		\]
\end{Theorem}

In this paper, we give an entirely geometric, and more conceptual proof of the classification of mutation-equivalence classes of $T$-polygons, and therefore also of Theorem~\ref{thm:main}.

As we explain below, the bijection between deformation families of smooth del Pezzo surfaces and mutation-equivalence classes of $T$-polygons should be viewed as an instance of mirror symmetry. In view of the open problem of classifying those algebraic varieties which are mirror to Fano varieties, our geometric proof has a significant advantage: while the classification of polygons up to mutation becomes an intractible combinatorial problem in dimension greater than~$2$, we expect many of the geometric methods employed in this article to generalize to higher dimension. In drawing inspiration from both Fano mirror symmetry~\cite{AkhtarCoatesCorti} and Gross--Hacking--Keel mirror symmetry \cite{GrossHackingKeelMS} this paper also contributes to the ongoing effort of reconciling these two flavors of mirror symmetry.

Our proof proceeds by taking a hint from mirror symmetry:
recall that the mirror to an $n$-dimensional Fano variety $X$ with an anticanonical divisor $E$ is expected to be a log Calabi--Yau surface $U$, with a regular function $W \colon U \rightarrow \CC$. The variety $U$ is usually called a Landau--Ginzburg model, and $W$ is called the superpotential. $U$ has many torus charts $j \colon (\CC^\times)^n \hookrightarrow U$; given
such a torus chart we obtain a Laurent polynomial $f$ by restricting the superpotential to the image of $j$. The toric variety $X_P$ associated to the Newton polygon $P$ of $f$ as defined above is expected to be a toric degeneration of $X$.
Moreover, the assertion that any two toric degenerations $X_P$ and $X_Q$ of $X$ can be connected by a chain of rational curves is closely related to the question whether the transition function $\varphi \colon (\CC^\times)^n \dashrightarrow (\CC^\times)^n$ between the two corresponding torus charts admits a factorization into cluster mutations. We recall here that a cluster mutation is a special kind of birational transformation of $(\CC^\times)^n$ which preserves the holomorphic volume form \[
\Omega=\left(\frac{1}{2\pi {\rm i}}\right)^n \frac{{\rm d}z_1}{z_1} \wedge \dots \wedge \frac{{\rm d}z_n}{z_n}.
\]
With this in mind, we give a brief summary of our geometric proof.
Given a $T$-polygon $P$, let~${\bigl(Y_P, \bar{D}\bigr)}$ be the toric surface associated to the normal fan of $P$, with $\bar{D}$ its toric boundary, note that $Y_P \setminus\bar{D}=(\CC^\times)^2$. By blowing up $Y_P$ in the base locus of a certain pencil, we construct a log Calabi--Yau pair $(Y, D)$ with an elliptic fibration $f \colon Y \rightarrow \PP^1$ such that $D$ is a singular fiber of type ${\rm I}_n$. We note that the complement $U=Y \setminus D$ is a log Calabi--Yau variety, and the restriction of $f$ to $U$ is a regular function.
We then use results of Friedman~\cite{Friedman} and the Torelli theorem of Gross--Hacking--Keel~\cite{GrossHackingKeelLooijengapairs} to show that the possible pairs $(Y, D)$ arising from this construction fall into $10$ isomorphism types, which precisely mirrors the 10 deformation families of smooth del Pezzo surfaces.
If two $T$-polygons $P$ and $Q$ give rise to the same pair $(Y, D)$ via this construction, we show that the induced birational map $\varphi$
\[
\begin{tikzcd}
	&(Y, D) \ar[rd] \ar[ld]&\\
	\bigl(Y_P, \bar{D}\bigr) \ar[rr, "\varphi", dashed] &&\bigl(Y_Q, \bar{D}\bigr)
\end{tikzcd}
\]
admits a factorization into cluster mutations, using Hacking--Keating~\cite[Proposition~3.27]{HackingKeating}. From there, it is not hard to show that $P$ is mutation-equivalent to $Q$.

The existence of the factorization is related to the classical theorem of Max Noether that a~plane birational map admits a factorization into Cremona transformations and automorphisms, with the important difference that we require the maps in the factorization to be volume preserving, i.e., to preserve $\Omega$.

Finally, it is worth noting that the bijection between deformation families of Fano varieties and mutation-equivalence classes of polygons is expected to continue to hold if $X$ is a Fano with log terminal singularities. For del Pezzo surfaces $X$ which are \emph{not} smooth (and admit a $qG$-degeneration to a toric del Pezzo surface), Corti \cite{Corti} has recently proved that there is a bijection between deformation families of $X$ and mutation-equivalence classes of the corresponding lattice polygons, thereby establishing Conjecture A of \cite{AkhtarCoatesCorti} for del Pezzo surfaces which are not smooth.

\section{Mutations of polygons and Laurent polynomials}\label{sec:combinatorics}
In this section, we give the necessary background on mutations and define maximally mutable Laurent polynomials.
	\subsection{Mutations of polygons}\label{sec:mutationsofpolygons}
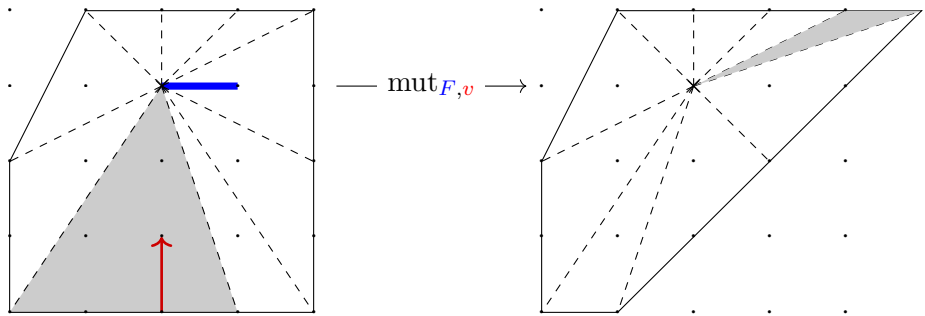
\begin{figure}[t]
		\centering
		\begin{tikzpicture}[scale=1]
			\begin{scope}[blend mode=multiply]
				\draw (-2,-1) -- (-2,-3) -- (2,-3) -- (2,1)--(-1,1)--(-2,-1);
				\draw[dashed] (0,0)--(2,1);
				\draw[dashed] (0,0)--(1,-3);
				\draw[dashed] (0,0)--(0,1);
				\draw[dashed] (0,0)--(1,1);
				\draw[dashed] (0,0)--(-1,1);
				\draw[dashed] (0,0)--(-2,-1);
				\draw[dashed] (0,0)--(-2,-3);
				\draw[dashed] (0,0)--(2,-3);
				\draw[dashed] (0,0)--(2,-1);
				\draw[line width=2.7pt, blue] (0,0)--(1,0);
				\draw[->, line width=1pt, red] (0,-3)--(0,-2);
				\fill[black!20] (0,0) -- (-2,-3) -- (1,-3) -- (0,0);
				\node at (-1,-2) {$\cdot$};
				\node at (0,-2) {$\cdot$};
				\node at (1,-2) {$\bullet$};
				\node at (-2,-2) {$\cdot$};
				\node at (-1,-1) {$\cdot$};
				\node at (0,-1) {$\cdot$};
				\node at (1,-1) {$\cdot$};
				\node at (2,-1) {$\cdot$};
				\node at (-2,-1) {$\cdot$};
				\node at (-2,0) {$\cdot$};
				\node at (-2,1) {$\cdot$};
				\node at (-1,0) {$\bullet$};
				\node at (0,0) {$\bullet$};
				\node at (1,0) {$\cdot$};
				\node at (2,0) {$\cdot$};
				\node at (-1,1) {$\cdot$};
				\node at (0,1) {$\cdot$};
				\node at (1,1) {$\cdot$};
				\node at (2,1) {$\cdot$};
				\node at (2,-2) {$\cdot$};
				\node at (0,-3) {$\cdot$};
				\node at (1,-3) {$\cdot$};
				\node at (2,-3) {$\cdot$};
				\node at (-2,-3) {$\cdot$};
				\node at (-1,-3) {$\cdot$};
			\end{scope}
			\begin{scope}[xshift=7cm]
				\draw (-1,1)--(-2,-1) -- (-2,-3) -- (-1,-3) -- (3,1)--(-1,1);
				\draw[dashed] (0,0)--(2,1);
				\draw[dashed] (0,0)--(0,1);
				\draw[dashed] (0,0)--(1,1);
				\draw[dashed] (0,0)--(-1,1);
				\draw[dashed] (0,0)--(1,-1);
				\draw[dashed] (0,0)--(3,1);
				\draw[dashed] (0,0)--(-1,-3);
				\draw[dashed] (0,0)--(-2,-1);
				\draw[dashed] (0,0)--(-2,-3);
				\fill[black!20] (0,0) -- (3,1) -- (2,1) -- (0,0);
				\node at (-1,-2) {$\bullet$};
				\node at (0,-2) {$\cdot$};
				\node at (1,-2) {$\cdot$};
				\node at (-2,-2) {$\cdot$};
				\node at (-1,-1) {$\cdot$};
				\node at (0,-1) {$\cdot$};
				\node at (1,-1) {$\cdot$};
				\node at (2,-1) {$\cdot$};
				\node at (-2,-1) {$\cdot$};
				\node at (-2,0) {$\cdot$};
				\node at (-2,1) {$\cdot$};
				\node at (-1,0) {$\bullet$};
				\node at (0,0) {$\bullet$}{};
				\node at (1,0) {$\cdot$};
				\node at (2,0) {$\cdot$};
				\node at (-1,1) {$\cdot$};
				\node at (0,1) {$\cdot$};
				\node at (1,1) {$\cdot$};
				\node at (2,1) {$\cdot$};
				\node at (2,-2) {$\cdot$};
				\node at (0,-3) {$\cdot$};
				\node at (1,-3) {$\cdot$};
				\node at (2,-3) {$\cdot$};
				\node at (-2,-3) {$\cdot$};
				\node at (-1,-3) {$\cdot$};
			\end{scope}
			\draw[->] (2.3,0) -- (4.8,0);
		\end{tikzpicture}
		\caption{Mutation of the polygon $P$ with respect to mutation data $\textcolor{red}{v}=(0,1)$, $\textcolor{blue}{F}=\Newt(1+x)$. The mutation contracts a grey $T$-cone on the left and adds a grey $T$-cone on the right. Residual lattice points in bold.}
		\label{fig:combinatorialmutation}
	\end{figure}

Let $M$ be a two-dimensional lattice with dual lattice $N$.
	\begin{Definition}\label{def:Fanopolygon}
		A Fano polygon is a full-dimensional lattice polygon $P \subset M_\RR$ such that $0$ is in the strict interior of $P$, and the vertices of $P$ are primitive lattice points.
	\end{Definition}
For any edge $E$ of $P$, the number of lattice points on $E$ minus one is called the \emph{lattice length}~$\ell_E$ of $E$. Let $u_E \in N$ be the primitive inward normal vector corresponding to $E \subset P$, then the positive integer $-\langle u_E, E \rangle$ is called the \emph{lattice height} $h_E$ of $E$.
$P$ is \emph{reflexive} if $h_E= 1$ for all edges $E$. Define $m_E$, $r_E$ to be the unique positive integers such that
\[
\ell_E=m_Eh_E+r_E, \qquad 0 \leq r_E<h_E.
\]
We call $r_E$ the \emph{residual length} of the edge $E$. Let $\sigma_E$ be the cone over the edge $E$.
If $r_E=0$, then $\sigma_E$ is called a $T$-cone. If $r_E=0$ and $m_E=1$ (or in other words $\ell_E=h_E$), $\sigma_E$ is called a primitive $T$-cone. If $\ell_E<h_E$, $\sigma_E$ is called a $R$-cone.
In general, we may (non-uniquely) subdivide $\sigma_E$ into $m_E$ primitive $T$-cones and zero or one $R$-cones, depending on whether $r_E$ is zero or nonzero.
Fixing such a subdivision, we say that a lattice point of $P$ is \emph{residual} if it is either the origin or interior to an $R$-cone. The number of residual points of $P$ is independent of the subdivision.
	\begin{Definition}
	A $T$-polygon is a Fano polygon such that every edge $E$ of $P$ satisfies $r_E=0$. Equivalently, the lattice length $\ell_E$ is divisible by the lattice height $h_E$.
\end{Definition}
As mentioned in the introduction, the classification of orbifold del Pezzo surfaces admitting a toric degeneration is conjecturally mirror to the classification of Fano polygons up to an appropriate equivalence relation. This equivalence relation is called mutation: while it is a bit technical to define, the idea behind it is rather simple, see Figure~\ref{fig:combinatorialmutation}.
	\begin{Definition}\label{mutationdef}
		Let $P \subset M$ be a Fano polygon and let $v \in N$ be a primitive vector. Choose a~line segment $F \subset v^{\perp} \subset M$ and write $P_d$
		for the slice of $P$ at height $d$ with respect to $v$. Suppose that for all $d<0$ we can decompose $P_d=R_d+(-d)F$ as a Minkowski sum for
		some line segment $R_d$ (where we allow $R_d =\varnothing$).
		Then we say that $P$ is mutable with respect to~$(v, F)$, and define the mutation of $P$ with respect to $(v, F)$ to be
		\[
		Q=\text{conv} \biggl( \bigcup_{d<0} R_d \cup \bigcup_{d\geq 0}(P_d+dF) \biggr).
		\]
		We call $F$ the factor of the mutation, and we say that two polygons $P$, $Q$ are {\it mutation equivalent} if there is a sequence of mutations
		of polygons starting with $P$ and ending with $Q$.
	\end{Definition}
	Keeping the notation of Definition~\ref{mutationdef}, suppose $F=kF'$ for some primitive line segment $F'$ and positive integer $k$ and suppose
	that $P$ is mutable with respect to $(v, F)$. Informally, $Q$ is obtained from $P$ by contracting $k$ primitive $T$-cones on one edge of $P$ and adding $k$
	primitive $T$-cones on the opposite edge of $P$.
	In particular, the condition for $P$ to be mutable with respect to $(v, F)$ means that there is an edge of $P$ perpendicular to $v$
	long enough to allow for the contracting of $d$ copies of $F$, where $d$ is the height of $d$ with respect to $v$ (see Figure~\ref{fig:combinatorialmutation} for an example with $k=1$ and $d=3$)
	By \cite[Proposition 3.6]{AkhtarKasprzyk}, any mutation of a $T$-polygon again produces a $T$-polygon, so that mutation defines an equivalence relation on the set of $T$-polygons.
	\subsection{Mutations of Laurent polynomials}\label{sec:mutationsofLaurent}
	Let now $v \in N$ and $f \in \CC[v^\perp] \subset \CC[M]$. Following \cite{GrossHackingKeelClusterAlgebras} and \cite{AkhtarCoatesCorti}, we define the automorphism~\smash{$x^{m} \!\mapsto x^{m}f^{\langle {m}, {v} \rangle}$}
	of the function field $\CC(M)$. This
 induces a birational map
${\varphi_{v,f} \colon T_N \dashrightarrow T_N}$
	which we call an \emph{algebraic mutation}. We call $f$ the \emph{factor} of the mutation. We will often suppress~$v$ and $f$ from notation.
	\begin{Definition}\label{def:AlgebraicMutation}
		Given a Laurent polynomial $g \in \CC[M]$, we say that $g$ is {\it mutable} with respect to an algebraic mutation $\varphi$ if $\varphi^*(g)\in \CC[M]$, i.e., $\varphi^*(g)$ is again a Laurent polynomial, and call~$\varphi^*(g)$ a \emph{mutation} of $g$.

		Given $g, g' \in \CC[M]$, we say that $g$ and $g'$ are \emph{mutation equivalent} if there exist algebraic mutations $\varphi_i$ for $1 \leq i \leq n$ and Laurent polynomials $g_i \in \CC[M]$ for $0 \leq i \leq n$ such that~${g_0=g}$, ${g_n=g'}$ and $\varphi_i^*g_{i-1}=g_{i}$ for all $i$.
	\end{Definition}

	Let us interpret mutability more concretely. Fix mutation data $v$ and $f$ as before. Taking the inner product with $v$ gives a $\ZZ$-grading of $M$ by height, so we may write $g=\sum_{d=-h}^m g_d$ where~$g_d$ is the sum of the monomials of $g$ at height $d$.
	Extend $v$ to a basis $e_1=v$, $e_2$ for $N$ and set $x=x^{e_2^*}$ and $y=x^{e_1^*}$. The factor $f$ is then a Laurent polynomial in $x$, and $\varphi_f^*(g_d)=g_df^{d}$. It follows that $g \in \CC[M]$ is mutable with respect to $\varphi_f$ if and only if $g_{-d}$ is divisible by \smash{$f^{d}$} for~${d>0}$.
	If $f$ is a monomial, then every Laurent polynomial is mutable with respect to $\varphi_f$, and we call such mutations \emph{trivial}. If a factor is of the form $(\lambda + x^u)$ for some $\lambda \in \CC^\times$ and {primitive} $u \in M$, we call the mutation \emph{standard}. It is clear that any factor is a product of standard and trivial factors.
	It follows easily from the definitions that a mutation of Laurent polynomials induces a mutation of their Newton polygons.
	However, the converse is false, it is not true that every mutation of $\Newt(g)$ is induced by a mutation of $g$ as the following example shows:
		\begin{figure}
		\centering
		\begin{tikzpicture}[scale=1]
				\begin{scope}[blend mode=multiply]
				\draw (-2,-1) -- (-1,-2) -- (1,-2) -- (2,1)--(-2,1)--(-2,-1);
				\draw(-1,-2) node[below]{$a$};
				\draw(0,-2) node[below]{$b$};
				\draw(1,-2) node[below]{$c$};
				\node(-1,-2) {$\cdot$};
				\node(0,-2) {$\cdot$};
				\node(1,-2) {$\cdot$};
				\node at (-2,-2) {$\cdot$};
				\node at (-1,-1) {$\cdot$};
				\node at (-2,-1) {$\cdot$};
				\node at (-2,0) {$\cdot$};
				\node at (-2,1) {$\cdot$};
				\node at (0,-1) {$\cdot$};
				\node at (1,-1) {$\cdot$};
				\node at (2,-1) {$\cdot$};
				\node at (-1,0) {$\cdot$};
				\node at (0,0) {$\times$}{};
				\node at (1,0) {$\cdot$};
				\node at (2,0) {$\cdot$};
				\node at (-1,1) {$\cdot$};
				\node at (0,1) {$\cdot$};
				\node at (1,1) {$\cdot$};
				\node at (2,1) {$\cdot$};
				\node at (2,-2) {$\cdot$};
				\fill[black!20] (0,0) -- (-1,-2) -- (1,-2) -- (0,0);
			\end{scope}
		\end{tikzpicture}
	\caption{}
	\label{fig:NotInduced}
	\end{figure}
	\begin{Example}\label{ex:NotInduced}
		Consider the Fano polygon $P$ in Figure \ref{fig:NotInduced}. $P$ is mutable with respect to $v=(0,1)$ and $F=\Newt(1+x)$.
		We have that
$
		g_{-2}=y^{-2}\bigl(a+bx+cx^2\bigr)x^{-1}
$
		for some constants $a$, $b$, $c$.
		In order for the mutation $\text{mut}_{F, v}$ to be induced by an algebraic mutation of $g$, the Laurent polynomial $g$ would have to be mutable with respect to a standard factor $f=\lambda + x$. This is only possible if $f^2$ divides $g_{-2}$ which happens if and only if $b^2=4ac$.
	\end{Example}

We say that a Laurent polynomial $g$ is supported on $P$ if	$\Newt(g) \subset P$.
It is natural to make the following definition.

	\begin{Definition}\label{def:TveitenClass}
Let $P$ be a Fano polygon. A Laurent polynomial $g$ supported on $P$ is of \emph{Tveiten class} if every mutation of $P$ is induced by an algebraic mutation $\varphi$ of $g$.
	\end{Definition}
We remark that the notion of \emph{maximally mutable Laurent polynomials} is now usually reserved for a more restricted class of polynomials (see below), so we opted to name the polynomials in Definition~\ref{def:TveitenClass} in view of their detailed study in Tveiten's~\cite{Tveiten} work on period integrals.

	Let us investigate the consequences of this definition, keeping the same notation as before. Fix an edge $E$ of $P$ with inner normal $v$ and write $\ell_E=mh+r$ as before. $P$ is mutable with respect to $(v, kF)$ for all $1 \leq k \leq m$, where $F \subset v^\perp$ be a primitive line segment.
	These mutations can only be induced by an algebraic mutation of $g$
	if there exists a polynomial $f \in \CC[x]$ with~${\Newt(f)=mF}$ such that for all $0 \leq d \leq h$, $g_{-d}$ is divisible by $f^d$ in $\CC[M]$. This is quite restrictive: up to a unit in $\CC[M]$ we may write
	$f=\prod_{i=1}^{m} (\lambda_i + x)$ with $\lambda_i \neq 0$, so that we have (again up to a unit)
	\[
	g_{-d}=\prod_{i=1}^{m} (\lambda_i + x)^{d} \cdot r_{-d},
	\]
	 where $r_{-d} \in \CC[x]$.

	 We see from this that a Laurent polynomial $g$ of Tveiten class is mutable with respect to~$m$ (not necessarily distinct) standard factors $(\lambda_i + x)$ along the edge $E$, one for each primitive $T$-cone on $E$.
	Since $\text{deg}(r_{-h})=r<h$, it follows that any Laurent polynomial $g$ can be mutable with respect to a maximum of $m$ standard factors along $E$, this motivates the term \emph{maximally mutable} used for Laurent polynomials of Tveiten class in \cite{Tveiten}. However, following the now standard terminology, we reserve this notion for those Laurent polynomials where all of the factors have $\lambda_i \equiv 1$.
	
	\begin{Definition}\label{def:MMLP}
		Let $P$ be a Fano polygon. A Laurent polynomial $g$ supported on $P$ is \emph{maximally mutable} if every mutation of $P$ is induced by an algebraic mutation $\varphi$ of $g$ and, moreover, the factor $f$ of $\varphi$ can always be taken to be $f=(1+x^u)^k$ for a primitive generator $u \in \CC [v^\perp ]$ and some $k\in \ZZ_{>0}$.
	\end{Definition}

If $g$ is maximally mutable, we see that for $0 \leq d \leq h$, we must have up to a unit that
	\[
	g_{-d}=\prod_{i=1}^m (1+x)^{d} \cdot r_{-d}=(1+x)^{dm} \cdot r_{-d},
	\]
	where $r_{-d} \in \CC[x]$. In particular, if $P$ is a $T$-polygon, then $g_{-h}=c(1+x)^{mh}$ for some $c \in \CC$.
		\begin{Definition}\label{def:normalized}
		A maximally mutable Laurent polynomial is \emph{normalized} if the coefficient of $g$ at each vertex of $\Newt(g)$ is $1$, and the constant term of $g$ is $0$.
	\end{Definition}
We will need the following result, which is very similar to \cite[Proposition 3.7]{Coates2021}, except that we are also interested in non-normalized Laurent polynomials.
\begin{Theorem}\label{thm:MMLP}
	Let $P$ be a $T$-polygon. There is a unique normalized maximally mutable Laurent polynomial $g_P$ supported on $P$. Moreover, the set of maximally mutable Laurent polynomial supported on $P$ is the two parameter family generated by $g_P$ and
	the constant Laurent polynomial~$1$.
\end{Theorem}
\begin{proof}
\cite[Proposition 3.7]{Coates2021} shows that $P$ supports a unique normalized maximally mutable Laurent polynomial $g_P$.

Suppose now that $g$ is any maximally mutable Laurent polynomial supported on $P$. If the coefficient of $g$ at any vertex is nonzero, we may scale $g$ by a scalar $\lambda$ to make this coefficient $1$. The mutability condition then implies that $\tfrac{1}{\lambda}g$ has binomial edge coefficients, and \cite[Proposition~3.7]{Coates2021} applies to show that $g=\lambda g_P+\mu$ for some scalar $\mu$.
If the coefficient of $g$ at a vertex is zero, then the mutability condition implies that the coefficient of $g$ along any boundary lattice point of $P$ is zero as well. The same argument as in \cite[Proposition 3.7]{Coates2021} then shows that $P$ can only have a nonzero coefficient at the origin, showing that $g=\mu$ for some scalar $\mu$. This completes the proof.
\end{proof}

	\section{The geometry of maximally mutable Laurent polynomials}\label{sec:geometry}
In this section, we show that Laurent polynomials of Tveiten class and maximally mutable Laurent polynomials can naturally be identified with the global sections of a certain line bundle. While well known to experts, this does not seem to be in the literature.
We then construct the log Calabi--Yau pair $(Y, D)$ associated to a $T$-polygon $P$, and show that $Y$ has the structure of a rational elliptic surface with $D$ a fiber of type ${\rm I}_n$. We have chosen to state some of our results in more generality than necessary, in the hope that this will clarify the geometric viewpoint on maximally mutable Laurent polynomials. We assume the reader is familiar with the basics of toric geometry.

	Throughout, let $(Y_P, D_P)$ be the polarized toric surface associated to the normal fan $\Sigma_P$ of~$P$. The edges $E$ of $P$ correspond to the rays of $\Sigma_P$, which in turn correspond to the toric divisors~$D_E$ of $Y_{P}$. The distinguished ample divisor on $Y_P$ is defined as
$D_P=\sum_{E \subset P} h_ED_E$.
	We may resolve the singularities of ${Y}_P$ to obtain a smooth toric surface $\bar{Y}_P$, and we denote its toric boundary $\bar{D}$.
	There is a well-known isomorphism
	\begin{equation}\label{eq:sections}
		\Gamma(Y_P, \cO(D_P)) \cong \bigoplus_{m \in P \cap M} \CC x^m.
	\end{equation}
Any Laurent polynomial $g$ supported on $P$ defines a section of $\cO(D_P)$. Its vanishing locus $C$ is a compactification of the affine curve $g=0$ in the dense torus $(\CC^\times)^2 \subset Y_P$. If $\Newt(g)=P$, the curve $C$ does not contain any of the toric divisors and does not pass through any torus fixed point. In particular, it avoids all singular points of $Y_P$.
	It follows that the strict transform of~$C$ on $\bar{Y}_P$ is isomorphic to $C$.
	To analyze the geometry, we work locally: given an edge $E$ of $P$ with inner normal $v$, let $e_1=v$, $e_2$ be a basis for one of the two smooth cones $\sigma$ in the fan of~$\bar{Y}_P$ containing $\RR_{\geq 0}v$ and let $x=x^{e_2^*}$, $y=x^{e_1^*}$ be the corresponding basis of characters. Then $D_E$ has an open subset $U$ with coordinates $x$, $y$ in which $D_E$ has local equation $y=0$, and $x$ is a local coordinate on $D_E$ around the torus fixed point corresponding to the cone $\sigma$, giving an identification $D_E^{\rm int} \cong \CC^\times$. Any mutation factor $f \in \CC[v^\perp]$ can then be identified with a Laurent polynomial $f(x)$. In particular, a standard mutation factor is of the form $f(x)=\lambda+x$.
	
	\begin{Definition}\label{def:mutablecycle}
		Suppose that $g$ is a Laurent polynomial supported on $P$ which is mutable with factor $(\lambda +x)$, and let $A$ be the point on $D_E$ where $\lambda+x=0$, i.e., the point with local coordinates $(-\lambda, 0)$. Then we say that $g$ is mutable with respect to $A$.

		Similarly, if $g$ is mutable with factor $(\lambda+x)^m$ for some positive integer $m$, we say that $g$ is mutable with respect to $m A$.
		The set of all points with respect to which $g$ is mutable defines a~zero cycle $Z$ supported on the interior of the toric boundary of $Y_P$, called the mutable cycle of~$g$.
	\end{Definition}
	We see that a Laurent polynomial $g$ with $\Newt(g)=P$ is of Tveiten class if and only if the mutable cycle has exact degree $m_E$ along the edge $E$, the maximal possible. The polynomial~$g$ is maximally mutable if and only if in addition the mutable cycle $Z$ is supported on the points $-1 \in D_E^{\rm int} \cong \CC^\times$ for $E \subset P$.
We now introduce the language of Looijenga pairs, which will simplify our discussion.
Following~\cite{Friedman, GrossHackingKeelLooijengapairs}, a Looijenga pair is a smooth projective surface $Y$ together with a singular anticanonical divisor $D$ with at worst nodal singularities. The divisor~$D$ is either a~nodal curve, or a cycle of~$n$ rational curves. $Y$ is necessarily rational, so we have an isomorphism $\Pic(Y) \cong H^2(Y, \ZZ)$.
If $Y$ is a toric surface with $D=Y \setminus (\CC^\times)^2$ its toric boundary, then $(Y,D)$ is called a {\it toric pair}.
Given a Looijenga pair $(Y, D)$, there are two elementary operations to produce another Looijenga pair:
\begin{itemize}\itemsep=0pt
	\item
	Let $p \colon Y' \rightarrow Y$ be the blowup of $Y$ at a smooth point of $D$. Denoting by $D'$ the strict transform of $D$, the pair $\bigl(Y', D'\bigr)$ is again a Looijenga pair. The map $p$ is called an \emph{interior blowup}.
	\item
	Let $p \colon Y' \rightarrow Y$ the blowup of a node of $D$. Denoting by $D'$ the reduced inverse image of~$D$, the pair $\bigl(Y', D'\bigr)$ is again a Looijenga pair. The map $p$ is called a \emph{corner blowup}.
\end{itemize}
	In the literature, corner blowups are often called toric blowups. To avoid confusion, we reserve this term for the blowup of a toric surface along a torus fixed point (which is a special case of a~corner blowup).
	We say that $(Y, D)$ is positive definite (negative definite, semi-definite,~\dots) if the intersection matrix of the components of $D$ is positive definite (negative definite, semi-definite,~\dots). A matrix is strictly negative semi-definite if it is negative semi-definite but not negative definite.
We will often identify a point $A \in D^{\rm int}$ with the corresponding point $E \cap \tilde{D}$ on the strict transform $\tilde{D}$ (where $E$ is the exceptional divisor of an interior blowup). In particular, the blowup of $m A$ is understood as $m$ iterated blowups at $A$.
We now aim to give a~geometric interpretation of maximal mutability.

 Let $P$ be a Fano polygon, and let $C$ be the vanishing locus of a section of $\cO(D_P)$ on the smooth toric surface $\bigl(\bar{Y}_P, \bar{D}\bigr)$, and let $A$ be a point in the intersection $\bar{D}^{\rm int} \cap C$. Let $p \colon (Y, D) \rightarrow \bigl(\bar{Y}_P, \bar{D}\bigr)$ be the blow up $kA$, and denote the exceptional classes $E_i$, for $1 \leq i \leq k$.
\begin{Lemma}\label{lem:mutabilitygeometric}
	The curve $C$ is the vanishing locus of a section of the line bundle
$
	p^*\cO(D_P)-\sum_{i=1}^k h E_i
$
	 on $Y$ if and only if in local coordinates $x$, $y$ centered at $A$, the curve $C$ has an equation of the form
	\begin{equation}\label{eq:23}
		\sum_{i=1}^hc_i y^ix^{k(h-i)}+\text{$($terms of degree $>kh)$}=0
	\end{equation}
for some constants $c_i$ and where $\deg(x)=1$, $\deg(y)=k$.
\end{Lemma}
\begin{proof}
	We may write the defining equation of $C$ as $\sum_{i,j \geq 0} a_{i,j} x^i y^j=0$. $C$ extends to a section of~${p^*D_P-hE_i}$ if and only if $A$ is a point of multiplicity at least $h$, i.e., $a_{i,j}=0$ for $i+j<h$.
	Locally, the blowup of $A$ is given by the map $\pi_1 \colon (u,v) \mapsto (u, uv)$, and the corresponding section~$C_1$ of $p^*D_P-hE_i$ has equation
$\sum_{i,j \geq 0} a_{i,j} u^{i+j-h}v^j=0$
	 in these coordinates. The point lying over~$A$ has coordinates $(u,v)=(0,0)$, so we see that $C_1$ extends to a section of $p^*D_P-hE_1-hE_2$ if and only if $C_1$ has a point of multiplicity at least $h$ at $(0,0)$, i.e., $a_{i, j}=0$ for $i+2j<2h$.
	Continuing in a similar fashion, we see that $C$ extends to a section of~\smash{$p^*D_P-\sum_{i=1}^k h E_i$} if and only if $a_{i,j}=0$ for $i+kj<kh$.
	This happens if and only if $C$ is of the form \eqref{eq:23}.
\end{proof}

Let now $Z$ be a zero cycle, supported on $\bar{D}^{\rm int}$, and let $k_E$ be the degree of $Z$ along the edge~$\bar{D}_E$.
Define $\pi \colon \tilde{Y}_Z \rightarrow \bar{Y}_P$ to be the blowup of $Z$. The exceptional locus is a disjoint union of chains of $\PP^1$s. Each such chain is of the form $C_1+\dots +C_r$, where $C_1^2=-1$, $C_i^2=-2$ for $i>1$, and $r$ is the multiplicity of the corresponding point in $Z$. For $1 \leq i \leq r$, we define $E_i=C_1+\dots+C_i$, and call the $E_i$ the exceptional classes of the blowup. Note that $E_i^2=-1$ for all $i$. Define the line bundle
\[
L=\pi^*D_P-\sum_{E\subset P} \sum_{i=1}^{k_E} h_EE_i.
\]
	\begin{Theorem}\label{thm:MMLPassections}
There is a one-to-one correspondence between Laurent polynomials $g$ supported on $P$ whose mutable cycle contains $Z$, and global sections of $L$.
\end{Theorem}
\begin{proof}
	Let $A \in \Supp(Z) \cap D_E$. Recall that we may choose a local coordinate chart $U$ centered around a torus fixed point on $D_E$ such that $D_E$ has equation $y=0$ and $A=(-\lambda, 0)$ for some constant $\lambda$. By definition, the mutable cycle of $g$ contains $kP$ if and only if we can write
	\[
	g=	
	\sum_{i=-h_E}^{0}c_i{y^{i}}(\lambda+x)^{-ki}{r}_{i}(x)+\sum_{i=1}^m {y^{i+h_E}}{r}_i(x)
	\]
	for some constants $c_i$, and Laurent polynomials $r_i(x)$.
	Multiplying through by a monomial, we see that $g$ is mutable at $kA$ if and only if the corresponding section $C$ of $\cO(D_P)$ (under the isomorphism \eqref{eq:sections}) is locally defined by an equation of the form
	\[	
	\sum_{i=0}^{h_E}c_i{y^{i}}(\lambda+x)^{k(h_E-i)}{r}_{i}(x)+\sum_{i=1}^m {y^{i+h_E}}{r}_i(x)=0
	\]
	for some constants $c_i$ and polynomials $r_i(x)$.
	Applying the coordinate change $x \mapsto x-\lambda$, we see that $C$ is mutable at $kA$ if and only if it has a local equation of the form
	\[	
	\sum_{i=0}^{h_E}c_i{y^{i}}x^{k(h_E-i)}+\text{(terms of degree $>kh_E$)}=0.
	\]
	By Lemma~\ref{lem:mutabilitygeometric}, this happens if and only if $g$ extends to a section of $\pi^*D_P-\sum_{i=1}^k h_E E_i$. Repeating the argument for all $A \in \Supp(Z)$ shows that $g$ extends to a section of $L$ if and only if the mutable cycle of $g$ contains $Z$.
	\end{proof}

In particular, by taking $Z$ to be maximal -- in the sense that $Z$ has degree $m_E$ along $D_E$ -- we can view Laurent polynomials of Tveiten class with mutable cycle $Z$ as the space of global sections of a certain line bundle. The most important case for us if $Z$ is supported on the points~${[1:-1] \in D_E}$, which corresponds to the maximally mutable Laurent polynomials.
	\begin{Corollary}
		Let $P$ be a Fano polygon. Let $\bigl(\bar{Y}_P, \bar{D}\bigr)$ be the associated smooth toric surface. Let~\smash{$p \colon \bigl(\tilde{Y}_Z, D\bigr) \rightarrow \bigl(\bar{Y}_P, \bar{D})$} be the Looijenga pair obtained by $m_E$ blowups at the point ${[-1:1] \in D_E}$, for every edge $E \subset P$. Let $E_i$ denote the exceptional classes of the blowups. There is a~$1$-$1$ correspondence between maximally mutable Laurent polynomials $g$ supported on~$P$ and global sections of
$p^*D_P-\sum_{E \subset P}\sum_{i=1}^{m_E} h_E E_i$,
		where $E \subset P$ ranges over all edges of $P$.
	\end{Corollary}
\begin{proof}
	Immediate from Theorem~\ref{thm:MMLPassections}.
\end{proof}

We now study the pencil of sections $\Gamma_g \subset |D_P|$ generated by $g$ and the constant Laurent polynomial $1$. Note that the vanishing locus of the section corresponding to $1$ is precisely $D_P$.
\begin{Lemma}\label{lem:baseschemevsmutablecycle}
	Let $P$ be a $T$-polygon, and $g$ with $\Newt(g)=P$ a Laurent polynomial of Tveiten class. Then the base scheme of the pencil $\Gamma_g$ and the mutable cycle of $g$ are supported on the same points.
\end{Lemma}
\begin{proof}
	Let $D$ be the toric boundary of the singular toric surface $Y_P$. The zero scheme of ${1 \in \cO(D_P)}$ is the divisor $D_P$, hence the base scheme of $\Gamma_g$ is supported on the points $D \cap \{g=0\}$. As before, for a fixed edge $E$ of $P$, we have local coordinates $x$, $y$ in which $D_E$ has equation~${y=0}$ and
	\[
	g=\prod_{i=1}^{m} (\lambda_i + x)^{h} \cdot r(x) +\text{(monomials involving $y$)}
	\]
	for some polynomial $r(x)$, and since $P$ is a $T$-polygon, $r(x) \equiv 1$. It follows that the intersection~${D_E \cap \{g=0\}}$ consists of the points $(x, y)=(-\lambda_i, 0)$, which coincides with the support of the mutable cycle of $g$ on $D_E$.
\end{proof}

We note however that the mutable cycle does not coincide with the base scheme of $\Gamma_g$ as soon as $P$ has an edge at height greater than one.
\begin{Lemma}\label{lem:resolutionofsingularities}
Suppose that $P$ is a $T$-polygon, and that $g$ is a Laurent polynomial of Tveiten class with $\Newt(g)=P$. Suppose that a point $A$ appears in the mutable cycle of $g$ with multiplicity~$k$. Then
\begin{enumerate}\itemsep=0pt
	\item[$(1)$] the basepoint $A$ of $\Gamma_g$ can be resolved by a composition $p \colon \bigl(\tilde{Y}, \tilde{D}\bigr) \rightarrow \bigl(\bar{Y}_P, \bar{D}\bigr)$ of $k$ interior blowups at $A$ $($see Figure~{\rm\ref{fig:h-uplepoint}}$)$.
	\item[$(2)$] The strict transform $\tilde{D}_P$ is a member of $p_*^{-1}\Gamma_g$ and
	\[p_*^{-1}\Gamma_g \subset \left|p^*\cO(D_P)-\sum_{i=1}^{k} h_E E_i\right|.
	\]
\end{enumerate}
\end{Lemma}
\begin{proof}
$A$ lies in the interior of a toric divisor $D_E$ corresponding to an edge $E$, let $h=h_E$.
As before, we may choose local coordinates $x$, $y$ such that $D_E$ has equation $y=0$. Suppose $A$ corresponds to $\lambda \in \CC^\times \cong D_E$. By assumption, $g$ locally has an equation of the form
\[
\sum_{i=0}^h y^i(x-\lambda)^{k(h-i)}r_i(x)+O\bigl(y^{h+1}\bigr)=0.
\]
It follows that $\Gamma_g$ has a local equation
\[
s\left(\sum_{i=0}^h y^i(x-\lambda)^{k(h-i)}r_i(x)+O\bigl(y^{h+1}\bigr)\right)+ty^h=0
\]
for $[s:t]\in \PP^1$.
Locally, the blowup of $A$ is given by the map $p_1 \colon (u,v) \mapsto (u, v(u-\lambda))$. Therefore, the strict transform of the pencil is given by
 \[
s\left(\sum_{i=0}^h v^i(u-\lambda)^{(k-1)(h-i)}r_i(u)+O\bigl(v^{h+1}\bigr)\right)+tv^h=0.
 \]
We note that the only basepoint of $p_{1*}^{-1}\Gamma_g$ on the exceptional divisor $u=\lambda$ is the point ${(u,v)=(\lambda, 0)}$, the intersection of the exceptional divisor with the strict transform of $D_E$.
After~$k$ blowups at $(\lambda, 0)$, the strict transform $p_*^{-1}\Gamma_g$ is given by
 \[
 s\left(\sum_{i=0}^h v^ir_i(u)+O\bigl(v^{h+1}\bigr)\right)+tv^h=0.
 \]
 $P$ is a $T$-polygon and $g$ is of Tveiten class, therefore $r_0(u)$ is a product of factors of the form~${(u-\lambda_i)^h}$, with $\lambda_i \neq \lambda$. It follows that $r_0(\lambda) \neq 0$, and so $p_*^{-1}\Gamma_g$ has no basepoint on the exceptional divisor $u=\lambda$, which proves the first claim.
The strict transform of $D_P$ under~$p$ has equation $v^h=0$, which is the member of $p_{*}^{-1}\Gamma_g$ corresponding to $[s:t]=[0:1]$. Finally, the proof shows that each blowup removes a basepoint of multiplicity $h=h_E$ from $\Gamma_g$, so that%
	\[
	p_*^{-1}\Gamma_g \subset \left|p^*D_P-\sum_{i=1}^{k} h_E E_i\right|.\tag*{\qed}
\] \renewcommand{\qed}{}
\end{proof}

\tikzset{every picture/.style={line width=0.75pt}} 
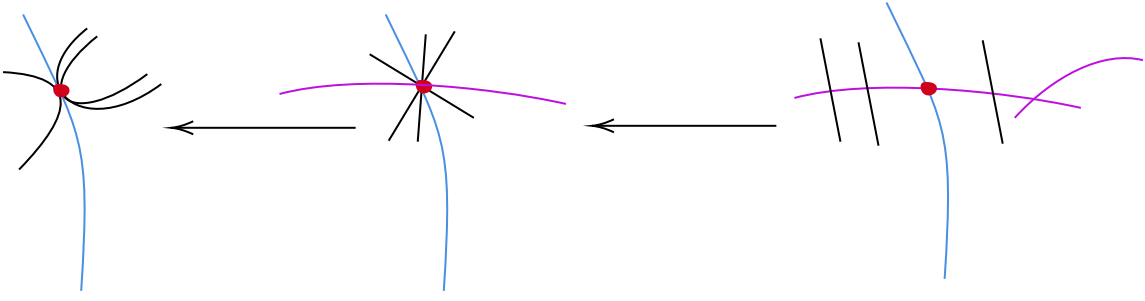
\begin{figure}[t]\centering
	\begin{tikzpicture}
		
		\draw [blue] (-2,0) -- (2,0);
		\draw (-2,-1.5) .. controls (-1,0.5) and (1,0.5).. (2,-1.5);
		\draw (-1.5,-1.5) .. controls (-1,0.5) and (1,0.5).. (1.5,-1.5);
		\draw (-2,1.5) .. controls (-1,-0.5) and (1,-0.5) .. (2,1.5);
		\filldraw [red] (0,0) circle (2pt);
		\draw [blue] (3.5,0) -- (7.5,0);
		\draw [green] (5.5,-1.5) -- (5.5,1.5);
		\draw (4.5, -1) -- (6.5, 1);
		\draw (5, -1) -- (6, 1);
		\draw (4.5, 1) -- (6.5, -1);
		\filldraw [red] (5.5,0) circle (2pt);
		\draw [blue] (9,0) -- (13,0);
		\draw [green] (11,-1.5) -- (11,1.5);
		\draw (10, -0.5) -- (12, -0.5);
		\draw (10, -1) -- (12, -1);
		\draw [green] (10, 1) -- (13, 1);
		\draw (10, 0.5) -- (12, 0.5);
		\filldraw [red] (11,0) circle (2pt);
		\draw [->] (3.25,0) --(2.25,0);
		\draw [->] (8.75,0) --(7.75,0);
	\end{tikzpicture}
		\caption{A triple point with two-fold tangency (along the blue divisor) and its resolution, exceptional divisors in green.}
\label{fig:h-uplepoint}
\end{figure}

We deduce the following theorem.
\begin{Theorem}
	Let $P$ be a $T$-polygon, and let $g$ a maximally mutable Laurent polynomial with~${\Newt(g)=P}$. Let $p \colon \tilde{Y}_Z \rightarrow \bar{Y}_P$ be the blowup of the mutable cycle $Z$ of $g$.
	\begin{enumerate}[$(a)$]\itemsep=0pt
	\item $p$ is a composition of interior blowups, and in particular $\bigl(\tilde{Y}_Z, \tilde{D}\bigr)$ is a Looijenga pair.
	\item The linear system $p_*^{-1}\Gamma_g$ defines a morphism $\pi \colon \tilde{Y}_Z \rightarrow \PP^1$ which is an elliptic fibration with connected fibers, and $\pi^*(\infty)=\tilde{D}_P$.
	\end{enumerate}
\end{Theorem}
\begin{proof}
	(a) is clear from the definition of $\tilde{Y}_Z$.

	For (b), applying Lemma~\ref{lem:resolutionofsingularities} repeatedly shows that
	\begin{equation}\label{eq:stricttransformofpencil}
	p_*^{-1}\Gamma_g \subset 	\left|p^*D_P-\sum_{E \subset P}\sum_{i=1}^{m_E} h_E E_i\right|.
	\end{equation}
	Lemma~\ref{lem:baseschemevsmutablecycle} shows that the only basepoints of $\Gamma_g$ occur at the support of the mutable cycle, and Lemma~\ref{lem:resolutionofsingularities} shows that blowing up the mutable cycle removes these basepoints, so that $p_*^{-1}\Gamma_g$ is base-point-free, and defines a morphism $\pi \colon \tilde{Y}_Z \rightarrow \PP^1$. By Bertini's theorem, the general member of $p_*^{-1}\Gamma_g$ is smooth, so the general fiber of $\pi$ is smooth. By Theorem~\ref{thm:MMLPassections}, sections of $L$ correspond to maximally mutable Laurent polynomials supported on $P$, and by Theorem~\ref{thm:MMLP}, the space of maximally mutable Laurent polynomials is two-dimensional. It follows that the inclusion \eqref{eq:stricttransformofpencil} is an equality, and $p_*^{-1}\Gamma_g$ is a complete linear system. Therefore, the morphism~$\pi$ is equal to its Stein factorization, so that $\pi$ has connected fibers. The fact that the general fiber has genus~$1$ follows from Riemann--Roch, see \cite[Theorem 3.5]{Tveiten}. Finally, Lemma~\ref{lem:resolutionofsingularities} shows that~$\tilde{D}_P$ is a~member of $p_*^{-1}\Gamma_g$, so that $\pi^*(\infty)=\tilde{D}_P$ after choosing suitable coordinates on~$\PP^1$.
\end{proof}

	In what follows, we will often denote the surface $\tilde{Y}_Z$ just constructed by $\tilde{Y}_g$ to emphasize the dependence of the morphism $\tilde{Y}_g \rightarrow \PP^1$ on the particular Laurent polynomial $g$. $\tilde{Y}_g$ fits in a~diagram
\[
\begin{tikzcd}
	(\CC^\times)^2 \ar[r] \ar[d, "g"]&\tilde{Y}_g\ar[d, "\pi"]\\
	\CC \ar[r]& \PP^1
\end{tikzcd}
\]
such that the fibre of $\pi$ over $[s:t] \in \PP^1$ is a compactification of the curve $g=\tfrac{s}{t}$. This fibration might not be relatively minimal, but we can contract all $(-1)$-curves contained in fibres to obtain a relatively minimal fibration $Y_g \rightarrow \PP^1$. We have the following.
\begin{Lemma}\label{lem:Every-1isToricBlowdown}
	Any $(-1)$-curve contained in a fibre of $\tilde{Y}_g$ is a component of $\pi^{*}(\infty)$.
\end{Lemma}
\begin{proof}
	Let \smash{$\tilde{D} \in |{-}K_{\tilde{Y}_g}|$} denote the strict transform of the toric boundary $\bar{D}$ under $\tilde{Y}_g \rightarrow \bar{Y}_P$. If $C$ were a $(-1)$-curve contained in a fibre different from $\pi^*(\infty)$, we would have $\tilde{D} \cdot C=0$, since $\tilde{D}$ is the underlying reduced curve of $\pi^*(\infty)$. However, adjunction gives $\tilde{D} \cdot C=1$, a~contradiction.
\end{proof}

It follows that passing from $\tilde{Y}_g$ to $Y_g$ amounts to a sequence of corner blowdowns $\bigl(\tilde{Y}_g, \tilde{D}\bigr) \rightarrow (Y_g, D)$.

Using Kodaira's classification of singular fibres of genus $1$ fibrations with connected fibers, (see, for example, \cite[Section~{\rm V}.7]{BarthPetersVandeVen}), we see that $D=mD'$ where $D'$ is either an irreducible nodal rational curve or a cycle of $n$ reduced rational $(-2)$-curves. Put differently, the intersection matrix of $D'$ must be strictly negative semi-definite.
In fact, we will show in Corollary~\ref{prop:m=1} that~${1 \leq n \leq 9}$ and that $D$ cannot be a multiple fiber (i.e., $m=1$) if $g$ is maximally mutable. We can thus summarize our findings as follows:
	\begin{Definition}\label{def:mainconstruction}
	Let $g$ be a Laurent polynomial of Tveiten class with $P=\Newt(g)$ a $T$-polygon, let $Y_P$ be the toric variety defined by $P$, and let $\bar{Y}_P$ be the toric minimal resolution. Then the associated strictly negative semi-definite Looijenga pair $(Y_g, D)$ is constructed by
	\begin{itemize}\itemsep=0pt
		\item Blowing up the mutable cycle on $\bar{Y}_P$, yielding a toric model $\bigl(\tilde{Y}_g, \tilde{D}\bigr) \rightarrow \bigl(\bar{Y}_P, \bar{D}\bigr)$, with a~genus $1$ fibration $\pi \colon \tilde{Y}_g \rightarrow \PP^1$.
		\item Blowing down $(-1)$-curves contained in fibres of $\pi$, yielding a corner blowdown $\bigl(\tilde{Y}_g, \tilde{D}\bigr) \rightarrow (Y_g, D)$ with a relatively minimal genus $1$ fibration $\pi \colon {Y}_g \rightarrow \PP^1$, such that $D=\pi^*(\infty)$.
	\end{itemize}
\end{Definition}
\begin{Example}
	Consider the $T$-polygon $P$ shown in the left of Figure \ref{fig:ExampleofConstruction}. The unique normalized maximally mutable Laurent polynomial with Newton polygon $P$ is \smash{$g=y+\frac{1}{xy}+\frac{2}{y^2}+\frac{x}{y^3}$}.
	The~fan of the minimal resolution $\bar{Y}_P$ of the toric variety $Y_P$ is shown on the right of Figure \ref{fig:ExampleofConstruction}. The generic member of the pencil $\Gamma_g$ has one basepoint of multiplicity two along the edge of~$P$ of length two, and one basepoint of multiplicity one along the two other edges. Blowing up these three basepoints, we arrive at the toric surface $\tilde{Y}_g$, as shown in Figure~\ref{fig:Y_f} (adapted from~\cite{ducat}). The strict transform of the toric boundary $\bar{D}$ has ten components, whose self-intersection numbers are shown in Figure~\ref{fig:Y_f}. Note that the component of self-intersection $(-1)$ appears with multiplicity~$2$ in $\pi^*(\infty)$.
	To obtain $Y_g$, we contract the $(-1)$ curve in $\pi^*(\infty)$. The fibre over $\infty$ is now a cycle of nine $(-2)$-curves. The elliptic surface $Y_g \rightarrow \PP^1$ is well known as the modular elliptic surface associated to the congruence subgroup $\Gamma_1(3)$.
\end{Example}

\begin{figure}[t]
\centering
		\begin{tikzpicture}
			\begin{scope}[scale=1]
				\draw (-1,-1) -- (1,-3) -- (0,1) -- (-1,-1);
				\node at (-1,-3) {$\cdot$};
				\node at (0,-3) {$\cdot$};
				\node at (1,-3) {$\cdot$};
				\node at (-1,-2) {$\cdot$};
				\node at (0,-2) {$\cdot$};
				\node at (1,-2) {$\cdot$};
				\node at (-1,-1) {$\cdot$};
				\node at (0,-1) {$\cdot$};
				\node at (1,-1) {$\cdot$};
				\node at (-1,0) {$\cdot$};
				\node (00) at (0,0) {$\times$}{};
				\node (10) at (1,0) {$\cdot$};
				\node at (-1,1) {$\cdot$};
				\node (01) at (0,1) {$\cdot$};
				\node at (1,1) {$\cdot$};
				\draw[dashed] (0,0)--(-1,-1);
				\draw[dashed] (0,0)--(1,-3);
				\draw[dashed] (0,0)--(0,1);
			\end{scope}
			\begin{scope}[scale=1, xshift=7cm]
				\draw[-stealth, red] (0,-1)--(-4,-2);
				\draw[-stealth, red] (0,-1)--(2,-2);
				\draw[-stealth] (0,-1)--(-3,-2);
				\draw[-stealth] (0,-1)--(-2,-2);
				\draw[-stealth] (0,-1)--(-1,-2);
				\draw[-stealth] (0,-1)--(0,-2);
				\draw[-stealth] (0,-1)--(1,-2);
				\draw[-stealth] (0,-1)--(1,-1);
				\draw[-stealth] (0,-1)--(-1,-1);
				\draw[-stealth, red] (0,-1)--(1,0);
				\node at (-4,-1) {$\cdot $};
				\node at (-3,-1) {$\cdot $};
				\node at (-2,-1) {$\cdot $};
				\node at (-1,-1) {$\cdot$};
				\node at (0,0) {$\cdot$};
				\node at (1,-1) {$\cdot$};
				\node at (2,-1) {$\cdot$};
				\node at (-4,0) {$\cdot $};
				\node at (-3,0) {$\cdot $};
				\node at (-2,0) {$\cdot$};
				\node at (-1,0) {$\cdot$};
				\node at (0,-1) {$\times$}{};
				\node at (1,0) {$\cdot$};
				\node at (2,0) {$\cdot$};
				\node at (-4,1) {$\cdot $};
				\node at (-3,1) {$\cdot $};
				\node at (-2,1) {$\cdot$};
				\node at (-1,1) {$\cdot$};
				\node at (0,1) {$\cdot$};
				\node at (1,1) {$\cdot$};
				\node at (2,1) {$\cdot$};
				\node at (-4,-2) {$\cdot $};
				\node at (-3,-2) {$\cdot $};
				\node at (-2,-2) {$\cdot$};
				\node at (-1,-2) {$\cdot$};
				\node at (0,-2) {$\cdot$};
				\node at (1,-2) {$\cdot$};
				\node at (2,-2) {$\cdot$};
				\node at (-4,-3) {$\cdot $};
				\node at (-3,-3) {$\cdot $};
				\node at (-2,-3) {$\cdot$};
				\node at (-1,-3) {$\cdot$};
				\node at (1,-3) {$\cdot$};
				\node at (2,-3) {$\cdot$};
			\end{scope}
		\end{tikzpicture}
		\caption{On the left a $T$-polygon $P$, divided into primitive $T$-cones. On the right the fan of the minimal resolution $\bar{Y}_P$, with the rays of the fan of $Y_P$ in red.}
		\label{fig:ExampleofConstruction}
\end{figure}
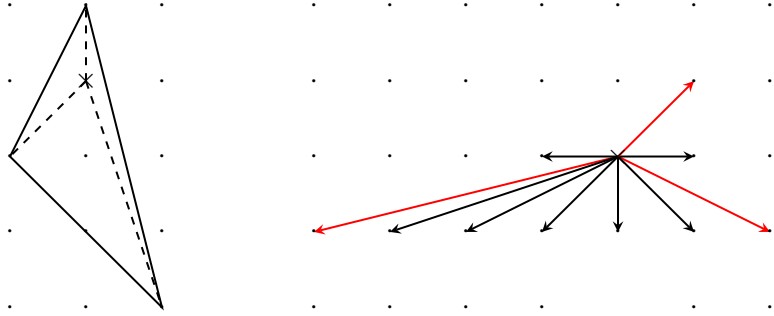
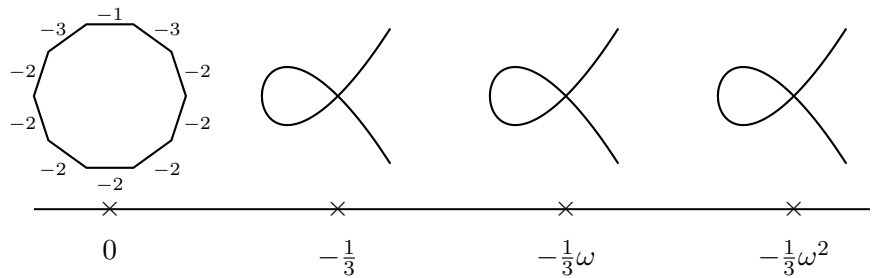
\begin{figure}[t]
	{\centering
		\begin{tikzpicture}
			\draw[thick] ({cos(0)},{sin(0)}) -- ({cos(36)},{sin(36)}) --({cos(72)},{sin(72)}) -- ({cos(108)},{sin(108)})--({cos(144)},{sin(144)}) -- ({cos(180)},{sin(180)}) --({cos(216)},{sin(216)}) -- ({cos(252)},{sin(252)}) --({cos(288)},{sin(288)}) -- ({cos(324)},{sin(324)}) --cycle;
			\draw[thick,scale=1,domain=0:1.3,smooth,variable=\x] plot ({3+\x^2-1},{\x^3-\x});
			\draw[thick,scale=1,domain=-1.3:0,smooth,variable=\x] plot ({3-1-\x^2},{\x^3-\x});
			\draw[thick,scale=1,domain=0:1.3,smooth,variable=\x] plot ({6+\x^2-1},{\x^3-\x});
			\draw[thick,scale=1,domain=-1.3:0,smooth,variable=\x] plot ({6-1-\x^2},{\x^3-\x});
			\draw[thick,scale=1,domain=0:1.3,smooth,variable=\x] plot ({9+\x^2-1},{\x^3-\x});
			\draw[thick,scale=1,domain=-1.3:0,smooth,variable=\x] plot ({9-1-\x^2},{\x^3-\x});
			\draw (-1,-1.5) -- (10,-1.5);
			\node at (3,-1.5) [label={below:$-\tfrac{1}{3}$}]{$\times$};
			\node at (6,-1.5) [label={below:$-\tfrac{1}{3}\omega$}]{$\times$};
			\node at (0,-1.5) [label={below:$0$}] {$\times$};
			\node at (9,-1.5) [label={below:$-\tfrac{1}{3}\omega^2$}] {$\times$};
			\node at (0,-0.8) [label={below:{\tiny $-2$}}] {};
			\node at (-0.75,-0.6) [label={below:{\tiny $-2$}}] {};
			\node at (0.75,-0.6) [label={below:{\tiny $-2$}}] {};
			\node at (1.15,0) [label={below:{\tiny $-2$}}] {};
			\node at (-1.15,0) [label={below:{\tiny $-2$}}] {};
			\node at (-1.15,0.7) [label={below:{\tiny $-2$}}] {};
			\node at (1.15,0.7) [label={below:{\tiny $-2$}}] {};
			\node at (0.75,1.25) [label={below:{\tiny $-3$}}] {};
			\node at (-0.75,1.25) [label={below:{\tiny $-3$}}] {};
			\node at (0,1.45) [label={below:{\tiny $-1$}}] {};
		\end{tikzpicture}
		\caption{The singular fibres of $\tilde{Y}_f \rightarrow \PP^1$.}
		\label{fig:Y_f}
	}
\end{figure}

\section[The classification of T-polygons]{The classification of $\boldsymbol{T}$-polygons}

	\subsection{Torelli for Looijenga pairs}\label{sec:Torelli for Looijenga pairs}
In this subsection, we use results of \cite{Friedman} and \cite{GrossHackingKeelLooijengapairs} to show that the pairs $(Y_g, D)$ constructed from a maximally mutable Laurent polynomials $g$ with $\Newt(g)$ a $T$-polygon fall into 10 isomorphism types.
	\begin{Definition}
		Let $(Y, D)$ be a Looijenga pair. Define the lattice
		\[
		\Lambda=\{L \in \Pic(Y) \mid L \cdot D_i=0 \; \text{for all}\; i \}.
		\]
	A cyclic ordering of the components of $D$ induces a canonical identification $\Pic^0(D) \cong \CC^\times$ (see~\cite[Lemma~2.1]{GrossHackingKeelLooijengapairs}).
		The map
		\[
		\phi_Y \colon\ \Lambda \rightarrow \Pic^0(D) \cong \CC^\times, \qquad L \mapsto L_{|D}
		\]
		is called the {\it period point} $\phi_Y \in \Hom(\Lambda, \CC^\times)$ of $(Y, D)$.
	\end{Definition}
	Let $\pi \colon \cY \rightarrow S$ be a flat morphism, from a smooth threefold $\cY$ to a smooth curve $S$. Suppose that $\cD$ is a relative anticanonical divisor with normal crossings on $\cY$ (i.e., $\cD$ restricts to a nodal anticanonical divisor on each fibre of $\pi$). We say that $(\cY, \cD)$ is a \emph{family of Looijenga pairs} if the family $\pi_{|\cD}$ is locally trivial on $S$. In particular, this implies that each anticanonical divisor $D_s$ has the same number of components.
	The inclusion $\cY_s \subset \cY$ of a fiber is a homotopy equivalence, so given a path between two points $s$, $t$ in $S$, we obtain an isometry
	\[
	\Pic(\cY_s) \cong H^2(\cY_s) \cong H^2(\cY) \cong H^2(\cY_t) \cong \Pic(\cY_t),
	\]
	which we will refer to as parallel transport.
	The main result of \cite{GrossHackingKeelLooijengapairs} is that the period point determines a Looijenga pair in a deformation family up to isomorphism.
	\begin{Theorem} \label{thm:torelli}
		Let $(Y,D)$ and $\bigl(Y',D'\bigr)$ be deformation-equivalent Looijenga pairs and
		suppose that $
		 \phi_{Y'} \circ \mu = \phi_{Y}$
		 under an isometry $\mu \colon \Pic(Y) \to \Pic(Y')$ induced by parallel transport.
		 Then there exists an isomorphism of pairs
		 $f\colon (Y,D) \to \bigl(Y',D'\bigr)$.
	\end{Theorem}
	\begin{Remark}
		In order to conclude $\mu=f^*$ in the statement above, one needs to add additional hypotheses, see \cite[Theorem 1.8]{GrossHackingKeelLooijengapairs} for details. For our purposes, the weaker conclusion of Theorem~\ref{thm:torelli} will be sufficient.
	\end{Remark}
	We now prove the main result of this section: if $f$ is maximally mutable, then the period point $\phi_{Y_f}$ of the surface $(Y_f, D)$ is equal to $1$. This follows easily from the results of \cite{GrossHackingKeelLooijengapairs}, but we first need a further definition.
	\begin{Definition}
		Given a Looijenga pair $(Y, D)$, a \emph{marking} of $D$ is a choice of point $p_i \in D_i^{\rm int}$ for all $i$. Given a marking of $D$, we define $\phi \in \Hom\bigl(\Pic(Y), \Pic^0(D)\bigr)$ by
		\[
		\phi(L)=(L_{\mid D}) \otimes \bigotimes_{i=1}^n \cO_D(-(L \cdot D_i)p_i).
		\]
		This is called the \emph{marked period point} of $(Y, D, p_i)$.
		Note that $\phi_{\mid \Lambda}=\phi_Y$ restricts to the period point of $Y$ as defined before.
	\end{Definition}
	
	\begin{Proposition}\label{periodpoint1}
		Let $P$ be a lattice polygon and $g$ be a maximally mutable Laurent polynomial with $\Newt P =g$. Then the Looijenga pair $(Y_g, D)$ associated to $g$ has period point $\phi_{Y_g}=1$, i.e., $\phi_{Y_g} \colon \Lambda \rightarrow \CC^\times$ is the constant function $1$.
	\end{Proposition}
	
	\begin{proof}
		Recall the diagram in Proposition~\ref{def:mainconstruction} summarising the construction of $(Y_g, D)$. We start with the toric pair $\bigl(\bar{Y}_P, \bar{D}\bigr)$ and then blow up the mutable cycle of $g$ to obtain a Looijenga pair~\smash{$\bigl(\tilde{Y}_g, \tilde{D}\bigr)$}.
		As before, a choice of orientation of $\bar{D}$ gives rise to a canonical identification~${D_i^{\rm int} \cong \CC^\times}$. Let $m_i$ correspond to $(-1)$ under this identification. Let $\tilde{\phi} \in \Hom\bigl(\tilde{Y}_g, \CC^\times\bigr)$ be the marked period point of \smash{$\bigl(\tilde{Y}_g, \tilde{D}, m_i\bigr)$} (where we identify $m_i$ with the corresponding point on the strict transform $\tilde{D}$ of $D$).
		\cite[Lemma 2.8]{GrossHackingKeelLooijengapairs} shows that the marked period point of~${\bigl(\bar{Y}_P, \bar{D}, m_i\bigr)}$ is~$1$, so \smash{$\tilde{\phi}_{\mid {\bar{Y}_P}}=1$}.
		$g$ is maximally mutable, so the mutable cycle of $g$ is supported on the points~$m_i$. It follows that in the construction of $\bigl(\tilde{Y}_g, \tilde{D}\bigr)$ we only blow up the points $m_i$ or points on a~strict transform of $D$ lying over $m_i$, so for any exceptional curve $E$ meeting $D_i$ we have that $\tilde{\phi}(E)=\cO_{\tilde{D}}(m_i) \otimes \cO_{\tilde{D}}(-m_i)=\cO_{\tilde{D}}$. Since the exceptional curves together with~$\Pic\bigl(\bar{Y}_P\bigr)$ generate \smash{$\Pic\bigl(\tilde{Y}_g\bigr)$} we conclude that \smash{$\tilde{\phi}=1$} and hence also that \smash{${\phi}_{\tilde{Y}_g}=\tilde{\phi}_{\mid \Lambda}=1$}.
		Finally, recall that in order to pass to $(Y_g, D)$, we perform a composition of toric blowdowns~\smash{$\pi \colon\bigl(\tilde{Y}_g,\tilde{D}\bigr) \rightarrow (Y_g, D)$}. This gives a canonical identification between the lattices $\Lambda_{\tilde{Y}_g}$ and $\Lambda_{Y_g}$ via $\pi^*$ and an isomorphism~\smash{$\pi^* \colon \Pic^0(D) \rightarrow \Pic^0\bigl(\tilde{D}\bigr)$}. The period points are then the same in the sense that \smash{$\pi^* \circ \phi_{Y_g}=\phi_{\tilde{Y}_g} \circ \pi^*$}, so we have that $\phi_{Y_g}=1$ as well.
	\end{proof}
	
	\begin{Remark}
		The statement that the period point of $(Y_g, {D})$ is $1$ is equivalent to the statement that the mixed Hodge structure on $H^2(Y_g \setminus D)$ is of Hodge--Tate type, see \cite[Proposition~3.12]{Friedman}.
	\end{Remark}

	Recall from Definition~\ref{def:mainconstruction} that $Y_g$ admits a relatively minimal genus $1$ fibration $\pi \colon Y_g \rightarrow \PP^1$.
	
	\begin{Corollary}\label{prop:m=1}
		Suppose that $g$ is maximally mutable and let $\pi \colon Y_g \rightarrow \PP^1$ be the associated genus~$1$ fibration. Then the fibre $\pi^*(\infty)$ is of Kodaira type ${\rm I}_n$ for $1 \leq n \leq 9$. In particular, $\pi^*(\infty)=D$, and $\pi$ admits a section.
	\end{Corollary}
	
	\begin{proof}
		Suppose first that $g$ is maximally mutable. From the discussion after Lemma~\ref{lem:Every-1isToricBlowdown}, we have that $\pi^*(\infty)=mD$ for some $m>0$, and therefore $D^2=0$ and $[D] \in \Lambda$.
		Applying the period point, we obtain
$
		\phi_{Y_g}([D])=\cO(D)_{|D}$.
		Since $g$ is maximally mutable, $\phi_{Y_g} \equiv 1$, so in particular~$\cO(D)_{|D}$ is trivial.
		However,~\cite[{\rm III}, Lemma 8.3]{BarthPetersVandeVen} shows that if $\pi^{-1}(\infty)=mD$ then~$\cO(D)_{|D}$ is torsion of exact order $m$. We conclude that $m=1$. Equivalently, all components of $\pi^*(\infty)$ must have multiplicity $1$, so $\pi^*(\infty)$ is of Kodaira type ${\rm I}_n$.
		In particular, $-K_{Y_g}$ is the class of a fibre $[F]$ and therefore any irreducible $(-1)$-curve $E$ (for example, the exceptional divisor of the last blowup) on $Y_g$ satisfies $E \cdot F=1$ by adjunction and so $E$ is a section of $Y_g$. Finally, since the Picard rank of $Y_g$ is $10$, we see that $n \leq 9$.
	\end{proof}

	We end this section with Friedman's classification of Looijenga pairs with negative semi-definite $D$. A sketch proof can be found in \cite[Theorems 9.15 and~9.16]{Friedman}, and is worked out in great detail in \cite{LutzThesis}.

	\begin{Theorem}[\cite{Friedman}]\label{thm:10deformationfamilies}
	There exist $10$ deformation families of strictly negative semi-definite Looijenga pairs $(Y, D)$.
\end{Theorem}

\begin{figure}[t]
	
\centering	\vspace{0.5em}
	\renewcommand{\arraystretch}{0.8}
	\begin{tabular}{cc}
		\hypertarget{poly:1}{}$1$&
		\hypertarget{poly:2}{}$2$\\
		\includegraphics[scale=1]{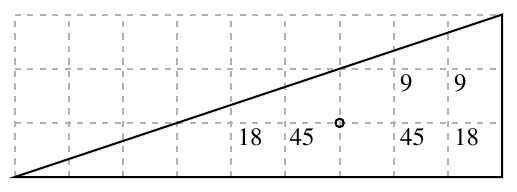}&
		\includegraphics[scale=0.24]{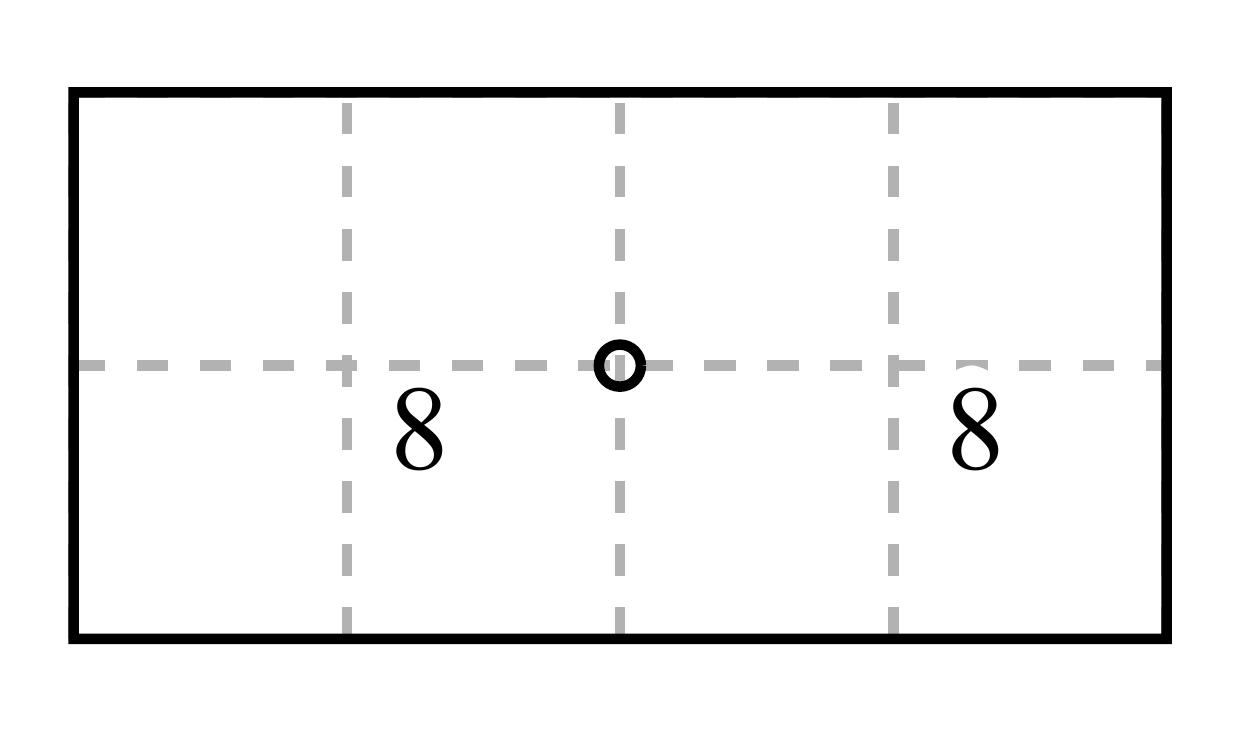}\\
	\end{tabular}
	\vgap
	
	\begin{tabular}{ccccc}
		\hypertarget{poly:3}{}$3$&
		\hypertarget{poly:4}{}$4$&
		\hypertarget{poly:5}{}$5$&
		\hypertarget{poly:6}{}$6$\\
		\includegraphics[scale=1]{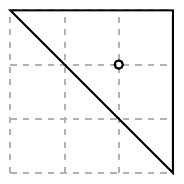}&
		\raisebox{13pt}{\includegraphics[scale=1]{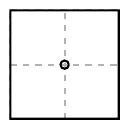}}&
		\raisebox{13pt}{\includegraphics[scale=1]{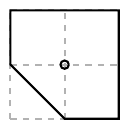}}&
		\raisebox{13pt}{\includegraphics[scale=1]{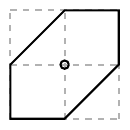}}\\
	\end{tabular}
	\vgap
	
	\begin{tabular}{cccccc}
		\hypertarget{poly:7}{}$7$&
		\hypertarget{poly:8}{}$8$&
		\hypertarget{poly:9}{}$8'$&
		\hypertarget{poly:10}{}$9$\\
		\raisebox{13pt}{\includegraphics[scale=1]{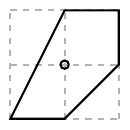}}&
		\raisebox{13pt}{\includegraphics[scale=1]{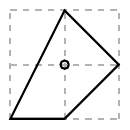}}&
		\raisebox{13pt}{\includegraphics[scale=1]{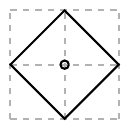}}&
		\raisebox{13pt}{\includegraphics[scale=1]{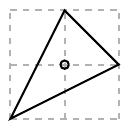}}\\
	\end{tabular}
\caption{Representatives for the ten mutation equivalence classes of $T$-polygons with their associated unique normalized maximally mutable Laurent polynomials, which are obtained by assigning binomial coefficients to lattice points on the edges of the polygon, $0$ to the origin, and the coefficients specified in the figure to the remaining lattice points. Picture taken from \cite{AkhtarCoatesCorti}, with a small correction in the coefficients of polygon $2$ (the labeled coefficients should both be $8$ instead of $4$).}	\label{T-polygonlist}
\end{figure}

	It is easy to write down representatives for each of the ten families, for example one can take the Looijenga pairs $(Y_{g_r}, D_{g_r})$ associated to the unique normalized maximally mutable Laurent polynomial $g_r$ with $\Newt(g_r)=P_r$, $1 \leq n \leq 10$ in Figure \ref{T-polygonlist}.

Given any $T$-polygon $P$, we have constructed a strictly negative semi-definite Looijenga pair~$(Y_g, D)$. By Theorem~\ref{thm:10deformationfamilies}, this Looijenga pair is deformation equivalent to one of the~$10$ reference pairs $(Y_{g_n}, D)$, and since both pairs have period point $1$, we actually have an isomorphism~${(Y_g, D) \cong (Y_{g_n}, D)}$.
	This gives rise to a diagram
\[
	\begin{tikzcd}
		&(Y_g, D) \cong (Y_{g_n}, D) \arrow[dr, "p"] \arrow[dl, swap, "q"] & \\
		\bigl(\bar{Y}_{P}, D\bigr) \ar[rr, dashed, "\varphi"]&&\, \bigl(\bar{Y}_{P_n}, D\bigr),
	\end{tikzcd}
	\]
	where $\varphi=p \circ q^{-1}$. In order to complete the classification of $T$-polygons, it remains to show that~$P$ is mutation-equivalent to $P_n$, which is the aim of the next subsection.

	\subsection{Volume-preserving birational maps}\label{sec:sarkisov}
In this subsection, we use the Sarkisov algorithm for volume-preserving maps \cite[Proposition~3.27]{HackingKeating} to complete our proof of the classification of $T$-polygons.
We start with a few preliminary results on volume-preserving birational maps.

 Let $(Y, D)$ be a Looijenga pair. Since $K_Y+D \sim 0$, there exists a nowhere vanishing volume form $\Omega$ on $Y \setminus D$ with simple poles along $D$, necessarily unique up to scaling.
	\begin{Definition}
		Let $(Y, D)$ and $\bigl(Y', D'\bigr)$ be Looijenga pairs, and let $\varphi \colon Y \dashrightarrow Y'$ be a birational map.
		We say that $\varphi$ is volume-preserving if there exists a resolution
		\[
		\begin{tikzcd}
			&Z \ar[dl, swap, "p"] \ar[dr, "p'"]&\\
			Y \ar[rr, dashed, "\varphi"]&&Y'
		\end{tikzcd}
	\]
		such that $p^*\Omega = \lambda p'^*\Omega'$ for some
		$\lambda \in \CC^\times$ where $\Omega$ is a holomorphic volume form on $Y \setminus D$ with simple poles along $D$, and similarly for $\Omega'$.
	\end{Definition}
We refer the reader to \cite[Remark 1.7]{CortiKaloghiros} for equivalent characterizations of this notion.
Volume-preserving morphisms have the following easy description:
\begin{Lemma}\label{lem:volume-preserving}
	Let $(Y, D)$ and $\bigl(Y', D'\bigr)$ be Looijenga pairs. A birational morphism $p \colon Y \rightarrow Y'$ is volume preserving if and only if it is a composition of corner blowups, interior blowups, and volume-preserving isomorphisms.
\end{Lemma}
\begin{proof}
	If $p$ is an interior or corner blowup, then one computes that ${p^*\bigl(K_{Y'}+D'\bigr)=K_Y+D}$, showing that $p$ is volume-preserving. For the converse, we can use \cite[Theorem 1.3.5]{Alberich} to factor~$p$ as
$
	p=p_n \circ \dots \circ p_1 \circ u$,
	where the $p_i$ are point blowups and $u$ is an isomorphism. We calculate that
	\[
	p^*\bigl(K_{Y'}+D'\bigr)=K_{Y}+u^*\tilde{D}-\sum_{i=1}^k u^*E_i,
	\]
	where the sum is over all exceptional divisors arising from blowups of points that are not on $D$ (or a strict transform of $D$), and $E_i$ is the pullback of the class of the corresponding exceptional divisor to $Y$.
	Since $p$ is volume preserving, we must have $u^*\tilde{D}=D$ and $k=0$, showing that $u$ is volume preserving, and that the $p_i$ are interior or corner blowups.
\end{proof}

\begin{Proposition}\label{thm:volume-preserving}
Let $(Y, D)$ and $\bigl(Y', D'\bigr)$ be Looijenga pairs, and $\varphi \colon Y \dashrightarrow Y'$ be a birational map. The following are equivalent
	\begin{enumerate}\itemsep=0pt
		\item[$(1)$] $\varphi$ is volume-preserving.
		\item[$(2)$] There exists a resolution
				\[
				\begin{tikzcd}
			&Z \ar[dl, swap, "p"] \ar[dr, "p'"]&\\
			Y \ar[rr, dashed, "\varphi"]&&\,Y',
		\end{tikzcd}
	\]
	where $p$ and $p'$ are volume-preserving birational morphisms

	\end{enumerate}
\end{Proposition}
\begin{proof}
	$(2) \implies (1)$: By Lemma~\ref{lem:volume-preserving}, $p$ and $p'$ are compositions of interior blowups, corner blowups and volume-preserving isomorphisms, so we have $p^*(K_Y+D)=K_{Z}+\tilde{D}$. It follows that $p^*(\Omega)$ is a nowhere vanishing holomorphic form on $Z \setminus \tilde{D}$ with simple poles along $\tilde{D}$. The same holds true for $p'^*\Omega$, so that $p'^*\Omega=\lambda p^*\Omega$ for some $\lambda \in \CC^\times$, i.e., $\varphi$ is volume-preserving.

	$(1) \implies (2)$: We have that $p^*(K_Y+D)=K_Z+\tilde{D}-\sum_i E_i$, where the sum is over all exceptional divisors arising from blowups of points that are not on $D$ (or a strict transform of~$D$), and $E_i$ is the pullback of the class of the corresponding exceptional divisor to $Z$. This means that there exists a holomorphic form on $Z \setminus \tilde{D}$ with simples poles along $\tilde{D}$ and simple zeros along the $E_i$.
	Similarly, we have $p'^*(K_{Y'}+D)=K_Z+\tilde{D'}-\sum_i F_i$. Since $\varphi$ is volume preserving, we must have $\tilde{D}=\tilde{D'}$ and $E_i=F_i$ up to reordering.

	 It follows that we may successively contract $(-1)$-curves on $Z$ to obtain a factorization
	\[
	\begin{tikzcd}
		&Z \ar[ddl, swap,"p"] \ar[ddr, "p'"] \ar[d]&\\
		&Z' \ar[dl, swap,"q"] \ar[dr, "q'"]&\\
		Y \ar[rr, dashed, "\varphi"]&&\,Y',
	\end{tikzcd}
	\]
	where $q^*(K_{Y}+D)=K_{Z'}+\tilde{D}$ and similarly for $q'$, meaning that $q$ and $q'$ are compositions of interior and corner blowups and hence volume-preserving.
\end{proof}

Recall that if $S \rightarrow C$ be a $\PP^1$-bundle over a curve and $p \in S$ a point, then the \emph{elementary transformation} at $p$ is the birational map $\alpha_p \colon S \dashrightarrow S'$ over $C$ given by blowing up $p$ and contracting the strict transform of the fiber through $p$.
\begin{Example}
If $S=\FF_k$ is a Hirzebruch surface and $k>0$, then $S'=\FF_{k+1}$ if $p$ lies on the negative section of $\FF_k$, and $S'=\FF_{k-1}$ otherwise.
If $S=\FF_0$, then $S'=\FF_1$.
\end{Example}
Considering $\FF_k$ as a Looijenga pair $(\FF_k, \partial \FF_k)$ (where $\partial \FF_k$ denotes its toric boundary), we see that the strict transform of $\partial \FF_k$ under an elementary transformation $\alpha_p \colon \FF_k \dashrightarrow \FF_{k \pm 1}$ is an anticanonical divisor if and only if $p$ lies on one of the torus-invariant sections of $\FF_k \rightarrow \PP^1$.
In particular, the birational map $\alpha_p \colon (\FF_k, \partial \FF_k) \dashrightarrow (\FF_{k \pm 1}, \partial \FF_{k \pm 1})$ is volume-preserving if and only if $p$ lies on one of the torus-invariant sections of $\FF_k \rightarrow \PP^1$, in which case we call $\alpha_p$ a \emph{mutation}.
	It is shown in \cite[Lemma 3.2]{GrossHackingKeelClusterAlgebras} that the restriction of the mutation $\alpha_p$ to the dense tori gives a map~${(\CC^\times)^2 \dashrightarrow (\CC^\times)^2}$ which is a mutation in the sense of Definition~\ref{def:AlgebraicMutation}.
	We will use the following theorem, which has appeared in \cite[Proposition 3.27]{HackingKeating}. A closely related version of this result was proved by \cite{Blanc}.
	\begin{Theorem}\label{thm:sarkisov}
		Let $(Y, D)$ be a Looijenga pair with two toric models
		\[
		\begin{tikzcd}
			&(Y, D) \ar[dl, swap, "p"] \ar[dr, "p'"] & \\
			\bigl(\bar{Y}, \bar{D}\bigr) \ar[rr, dashed, "\varphi"]&& \,\bigl(\bar{Y}', \bar{D}'\bigr).
		\end{tikzcd}
	\]
		Then $\varphi$ has a factorization
		\[
		\bigl(\bar{Y}, \bar{D}\bigr)=\bigl(\bar{Y}_0, \bar{D}_0\bigr) \xdashrightarrow{\varphi_1} \bigl(\bar{Y}_1, \bar{D}_1\bigr) \xdashrightarrow{\varphi_2} \cdots \xdashrightarrow{\varphi_n} \bigl(\bar{Y}_n, \bar{D}_n\bigr)=\bigl(\bar{Y}', \bar{D}'\bigr),
		\]
		 where each of the maps $\varphi_k$ is a toric blowup, toric blowdown, or a mutation.

		Moreover, let $p_k=\varphi_k \circ \dots\circ \varphi_1 \circ p$. Then $p_k \colon (Y, D) \dashrightarrow \bigl(\bar{Y}_k, \bar{D}_k\bigr)$ extends to a regular map~\smash{$\tilde{p}_k \colon \bigl(\widetilde{Y}, \widetilde{D}\bigr) \rightarrow \bigl(\bar{Y}_k, \bar{D}_k\bigr)$}
		 on some corner blowup \smash{$\bigl(\widetilde{Y}, \widetilde{D}\bigr)$} of $(Y, D)$.
	\end{Theorem}
	If $\bigl(\bar{Y}, \bar{D}\bigr)$ is a toric Looijenga pair and $p \colon (Y, D) \rightarrow \bigl(\bar{Y}, \bar{D}\bigr)$ is a composition of toric blowdowns and nontoric blowdowns then $p^{-1}\bigl((\CC^\times)^2\bigr)$ is a well-defined torus chart on $U=Y \setminus D$. Conversely, any torus chart arises in this way from a toric model. We conclude:
	
	\begin{Corollary}\label{cor:toruscharts}
		Any two torus charts $j, j' \colon (\CC^\times)^2 \hookrightarrow U$ on a Looijenga pair $(Y, D)$ with ${U=Y \setminus D}$ are related by a composition of algebraic mutations between torus charts on $U$.
	\end{Corollary}
\begin{proof}
	The two torus charts give rise to two toric models of a Looijenga pair $(Y, D)$ with $U=Y \setminus D$. Let $\varphi$ be the induced birational map between the two toric models.
The first part of Theorem~\ref{thm:sarkisov} gives a factorization $\varphi=\varphi_n \circ \dots \circ \varphi_1$. Restricting to the dense tori of the various toric models gives a commutative diagram
		\[
	\begin{tikzcd}
		&U && &\\
		(\CC^\times)^2 \ar[ur, "j"]\ar[r, dashed, "\varphi_n"]&\cdots\ar[r, dashed, "\varphi_1"]&\, (\CC^\times)^2 \ar[ul, swap, "j'"],&
	\end{tikzcd}
	\]
	where each $\varphi_i$ is now either an algebraic mutation (if the corresponding $\varphi_i \colon \bigl(\bar{Y}_{i-1}, \bar{D}_{i-1}\bigr) \xdashrightarrow{\varphi_i}\bigl(\bar{Y}_{i}, \bar{D}_{i}\bigr)$ was a mutation), or the identity (if the corresponding $\varphi_i$ was a toric blowup or toric blowdown).
	The second part of Theorem~\ref{thm:sarkisov} implies that the birational map
	\[
 \varphi_k \circ \dots \circ \varphi_1 \circ j' \colon\ (\CC^\times)^2 \hookrightarrow U
	\]
	 is a torus chart on $U$ for all $k$. It follows that the two torus charts $j$ and $j'$ are connected by a~composition of mutations between torus charts, as required.
\end{proof}

\begin{Proposition}\label{prop:volumepreserving=>mutation}
	Let $P$ and $Q$ be $T$-polygons, and let $\varphi \colon \bigl(\bar{Y}_P, \bar{D}\bigr) \dashrightarrow \bigl(\bar{Y}_Q, \bar{D}\bigr)$ be volume-preserving.
	Let $f$ and $g$ be Laurent polynomials of Tveiten class with $P=\Newt(f)$ and ${Q=\Newt(g)}$, and suppose that $\varphi^*f=g$. Then $f$ and $g$ are mutation-equivalent.
\end{Proposition}

\begin{proof}
By Proposition~\ref{thm:volume-preserving}, there exists a Looijenga pair $(Y, D)$ and a minimal resolution
	\[
	\begin{tikzcd}
		& (Y, D) \ar[dr, "q"] \ar[dl, swap, "p"]&\\
		\bigl(\bar{Y}_P, \bar{D}\bigr) \ar[rr, dashed, "\varphi"] &&\,\bigl(\bar{Y}_Q, \bar{D}'\bigr),
	\end{tikzcd}
	\]
	where $p$ and $q$ are volume-preserving morphisms. In particular, this means that $p$ only blows up points on $\bar{D}$ (or a strict transform thereof) and similarly for $q$. $f$ and $g$ define rational maps~${\phi_f \colon \bar{Y}_Q \dashrightarrow \PP^1}$ and $\phi_g \colon \bar{Y}_P \dashrightarrow \PP^1$. We will first show that $p$ and $q$ only blow up basepoints of $\phi_f$ and $\phi_g$. Indeed, $g$ takes the value $\infty$ on $\bar{D}$, so that if $p$ blows up a point which isn't a~basepoint of $\phi_g$, then $\phi_g$ would take the value $\infty$ on the corresponding exceptional divisor.
	In particular, there would exist a $(-1)$-curve $C$ on $Y$ with $g$ equal to $\infty$ along $C$.

	$C$ cannot be $q$-exceptional: otherwise, we can factor $q=q_1 \circ \dots \circ q_n$ as a composition of contractions of $(-1)$-curves $E_n, \dots, E_1$ with $C=E_i$ for some $1 \leq i \leq n$. Since $C$ is a $(-1)$-curve, any $E_j$ meeting $C$ must have $j<i$. It follows that the union of $E_n, \dots, E_{i+1}$ is disjoint from $C$, and therefore $q$ factors through the contraction of $C$. However, this contradicts the minimality of the resolution $(Y, D)$.

	Since $C$ is not $q$-exceptional, the intersection $q(C_n) \cap (\CC^\times)^2$ is nonempty, and since $f \circ \varphi=g$, we see that $f \circ q=g \circ p$ and therefore $f(q(C_n))=\infty$, a contradiction, since $f$ is a Laurent polynomial. We conclude that $p$ only blows up basepoints of $\phi_g$ and similarly $q$ only blows up basepoints of $\phi_f$.

	It follows that $\phi_f$ and $\phi_g$ extend to morphisms $\phi_f$ and $\phi_g \colon Y \setminus Z \rightarrow \PP^1$ (where $Z$ is the union of the remaining basepoints of the pencils defined by $f$ and $g$) which agree on a dense open subset (since $f \circ \varphi=g$) and are therefore equal. This implies that the remaining basepoints of~$\phi_f$,~$\phi_g$ are the same, and therefore after blowing them up, we have an isomorphism of elliptic surfaces $\tilde{Y}_f \cong \tilde{Y}_g$.

By Corollary~\ref{cor:toruscharts}, we obtain a factorization $\varphi=\varphi_n \circ \dots \circ \varphi_1$, which fits in a diagram
	\[
	\begin{tikzcd}
		&U \ar[rr, "W"] && \CC &\\
		(\CC^\times)^2 \ar[ur, "j"]\ar[r, dashed, "\varphi_n"]&\cdots\ar[r, dashed, "\varphi_1"]&\, (\CC^\times)^2 \ar[ul, swap, "j'"],&
	\end{tikzcd}
	\]
	where $U=\tilde{Y}_f \setminus D=\tilde{Y}_g \setminus D$, and $W$ is the restriction of the elliptic fibration $\phi_f=\phi_g$ to $U$. Note that Corollary~\ref{cor:toruscharts} shows that each $j_k=
	\varphi_k \circ \cdots \circ \varphi_1 \circ j'$ is a torus charts on $U$.
	By construction, we have $g=j^*W$ and $f=j'^*W$, and the pullback $f_k:=(\varphi_{1} \circ \dots \circ \varphi_k)^*f$ is similarly given by restricting $W$ via the torus chart $j_k=
	\varphi_k \circ \dots \circ \varphi_1 \circ j'$.
	In particular, $f_k$ is a regular function on~$(\CC^\times)^2$, so it is a Laurent polynomial, for all $k$, so by definition, each~$f_k$ is mutable with respect to $\varphi_{k+1}$. It follows that there is a sequence of algebraic mutations mapping ${g=f_0 \mapsto f_1 \mapsto \dots \mapsto f_n=g}$, which induces a sequence of mutations ${Q=\Newt(g) \rightarrow \dots \rightarrow \Newt(f)=P}$. This proves that~$P$ is mutation-equivalent to $Q$, as required.
\end{proof}

We can now complete the classification of $T$-polygons
\begin{Theorem}\label{cor:10mutationclasses}
	There are ten mutation equivalence classes of $T$-polygons.
\end{Theorem}
\begin{proof}
	Let $P$ be any $T$-polygon, and let $g$ be a maximally mutable Laurent polynomial with $\Newt(g)=P$. The pair $(Y_g, D)$ is strictly negative semi-definite by Definition~\ref{def:mainconstruction}, so by Theorem~\ref{thm:10deformationfamilies}, the Looijenga pair $(Y_g, D)$ must be deformation equivalent to one of the pairs $(Y_{g_n}, D_{g_n})$, where $g_n$ are as in Figure~\ref{T-polygonlist}. Since both pairs have period point $1$, the pairs must be isomorphic by the Torelli Theorem~\ref{thm:torelli}, so we have a diagram
	\[
	\begin{tikzcd}
		&Y_g=Y_{g_n} \arrow[dr] \arrow[dl] & \\
		(\CC^\times)^2 \subset \bar{Y}_{P} \ar[rr, dashed, "\varphi"]&&\,\bar{Y}_{P_n} \supset (\CC^\times)^2,
	\end{tikzcd}
	\]
	where the vertical maps are toric models, and $\varphi$ is volume-preserving.
	 It follows that the induced birational map $\varphi \colon (\CC^\times)^2 \dashrightarrow (\CC^\times)^2$ satisfies $\varphi^*\Omega=\lambda \Omega$ for some $\lambda \in \CC^\times$, where \smash{$\Omega=\bigl(\frac{1}{2\pi {\rm i}}\bigr)^2\frac{{\rm d}x \wedge {\rm d}y}{xy}$}. By pairing $\Omega$ with the integral generator $\{|x|=1, |y|=1\} \in H_2\bigl((\CC^\times)^2, \ZZ\bigr)$ and using the change of variable formula, we see that $\lambda=\pm 1$. Composing with the automorphism of $(\CC^\times)^2$ given by~$(x,y) \mapsto \bigl(x, \frac{1}{y}\bigr)$, we may assume that $\lambda=1$.

	Let $p$ and $q$ be the elliptic fibrations on $Y_g=Y_{g_n}$ induced by $g$ and $g_n$, note that by Proposition~\ref{prop:m=1}, the elliptic surface $Y_g$ has a section $s$. The isomorphism $(Y_g, D) \cong (Y_{g_n}, D)$ of Looijenga pairs maps the fibre $p^{-1}(\infty)$ to $q^{-1}(\infty)$. It follows from \cite[Lemma~1.5]{ClemensKollarMori} that the isomorphism maps every fibre of $p$ to a fibre of $q$ and therefore $\alpha=q \circ s$ is an automorphism of $\PP^1$ making the diagram
	\[
	\begin{tikzcd}
		&\PP^1 \\
		Y \arrow[dr, swap, "p"] \arrow[ur, "q"] &\\
		& \PP^1 \arrow[uu, "\alpha"]
	\end{tikzcd}
	\]
	commute.
	Since $\alpha$ fixes $\infty$, $\alpha$ must be of the form $z \mapsto az+b$ for $a,b \in \CC$ and $a \neq 0$.
	It follows that $\varphi^*{g_n}=ag+b$. Set $g':=ag+b$, then the Laurent polynomials $g_n$ and $g'$ satisfy the assumptions of Proposition~\ref{prop:volumepreserving=>mutation} and we conclude that $g_n$ and $g'$ are mutation equivalent. Therefore, $P=\Newt(g')$ is mutation-equivalent to $P_n=\Newt(g_n)$.
	It follows that any $T$-polygon $P$ is mutation-equivalent to one of the ten $T$-polygons in Figure~\ref{T-polygonlist}, which completes the proof.
\end{proof}

\subsection*{Acknowledgements}
I would like to thank Alessio Corti for suggesting the topic of this paper to me and for feedback on an earlier draft.
I would also like to thank Tom Coates, Robert Friedman, Paul Hacking, and Alan Thompson for helpful discussions.
I am grateful to the referees for their time and their helpful remarks which have improved my paper.
I was supported by the EPSRC Centre for Doctoral
Training in Geometry and Number Theory at the Interface, grant number EP/L015234/1.

\pdfbookmark[1]{References}{ref}
\LastPageEnding

\end{document}